\documentclass[12pt]{amsart}
\usepackage{latexsym,enumerate}
\usepackage{amsmath,amsthm,amsopn,amstext,amscd,amsfonts,amssymb,mathrsfs}
\usepackage{color}

\usepackage{color}
\usepackage{graphicx}
\usepackage{cite}
\numberwithin{equation}{section}

\newtheorem{theorem}{Theorem}
\newtheorem{lemma}{Lemma}
\newtheorem{corollary}{Corollary}
\newtheorem{proposition}{Proposition}
\newtheorem{remark}{Remark}

\numberwithin{theorem}{section}
\numberwithin{corollary}{section}
\numberwithin{lemma}{section}
\numberwithin{definition}{section}
\numberwithin{proposition}{section}
\numberwithin{remark}{section}

\newcommand{\R}{\mathbb R}

\newcommand{\medint}{-\kern  -,375cm\int}
\newcommand{\dint}{\displaystyle\int}

\begin{document}

\title[]{The isoperimetric problem for a class of non-radial weights and applications
 }

\author[A. Alvino]{A. Alvino$^1$}
\author[F. Brock]{F. Brock$^2$}
\author[F. Chiacchio]{F. Chiacchio$^1$}
\author[A. Mercaldo]{A. Mercaldo$^1$}
\author[M.R. Posteraro]{M.R. Posteraro$^1$}


\setcounter{footnote}{1}
\footnotetext{Universit\`a di Napoli Federico II, Dipartimento di Matematica e Applicazioni ``R. Caccioppoli'',
Complesso Monte S. Angelo, via Cintia, 80126 Napoli, Italy;\\
e-mail: {\tt angelo.alvino@unina.it, fchiacch@unina.it,  mercaldo@unina.it, posterar@unina.it}}

\setcounter{footnote}{2}
\footnotetext{South China University of Technology, International School of 
Advanced Materials, Wushan Campus 381, Wushan Road, 
Tianhe District, Guangzhou, P.R. China, 510641, email: fbrock@scut.edu.cn}

\begin{abstract} 
We study a class of isoperimetric problems on $\mathbb{R}^{N}_{+} $ where the densities of
the weighted volume and  weighted perimeter are given by two different non-radial functions of the type $|x|^k x_N^\alpha$.   
Our results  imply some sharp functional inequalities, like for instance, Caffarelli-Kohn-Nirenberg type inequalities. \medskip

\medskip

\noindent
{\sl Key words: isoperimetric inequality, weighted rearrangement, functional inequalities}  
\rm 
\\[0.1cm]
{\sl 2000 Mathematics Subject Classification:} 51M16, 46E35, 46E30, 35P15 
\rm 
\end{abstract}
\maketitle

\section{Introduction }
The last decades have seen an increasing  interest in the study of
 ``Manifolds with Density'', which is a manifold where both perimeter and
volume carry the same  weight.  To have an idea of the 
 possible applications of that subject one can consult, for instance \cite{Mo}, \cite{Mo1} 
 and the references therein.
  In particular, much attention has been devoted
 to find, for a given  manifold with density, its isoperimetric set (see, e.g., 
 \cite{BCMR}, \cite{BBCLT},  \cite{BCM, BCM2, BCM3, BMP, XR}, 
 \cite{CMV}, \cite{CJQW}, \cite{Cham}, \cite{DDNT}, \cite{Howe},
\cite{KZ}, \cite{MadernaSalsa},  \cite{Mo1}, \cite{Mo2}).
 On the other hand, many authors have studied isoperimetric problems
 when volume and perimeter carry two different weights.
 A remarkable example is obtained when
the manifold is $\mathbb R^{N}$ and the two weights are  two
different powers of the distance from the origin.
More precisely, given two real numbers $k$ and $l$, the problem is to find the set $G$ in $\mathbb R^{N}$
which minimizes the weighted perimeter $\displaystyle\int_{\partial G } |x|^k  \, {\mathcal H}_{N-1} (dx) $ 
once the weighted volume $\displaystyle\int_{ G } |x|^l  \, dx $ is prescribed.
Such a problem is far from being  artificial since its solution allows to compute, for instance,
 the best constants in the well-known Caffarelli-Kohn-Nirenberg inequalities as well as 
to establish the radiality of the corresponding minimizers. Several partial results have been obtained on such an issue (see, e.g., 
 \cite{ABCMP}, \cite{BBMP2},  \cite{C},   \cite{diGiosia_etal},  \cite{DHHT}, \cite{Howe},  \cite{Mo}) and a
complete solution is contained in in the recent paper 
(see \cite{diGiosia_etal}). There the authors find the full range of the 
parameters $k$ and $l$ for which the isoperimetric set is the ball centered at the origin.
The first step of their proof consists of reducing the problem
into a two-dimensional one by means of spherical symmetrization (also known as
foliated Schwarz symmetrization).
\\
Let 
 $\mathbb{R}^{N} _{+}  := \{ x \in \mathbb{R}^N :\, x_N >0\} $.  
The problem that we address here is the following:

Given  $k,l \in \mathbb{R}$, $\alpha > 0$,
\medskip

\noindent {\sl Minimize $\displaystyle\int_{\partial \Omega } |x|^k x_N^\alpha \, {\mathcal H}_{N-1} (dx) $ 
among all smooth sets 
$\Omega \subset \mathbb{R} ^{N}_{+}$ satisfying $\dint_{\Omega } |x|^lx_N^\alpha \, dx =1$.}
\medskip
 
Let $B_R$ denote the ball of $\mathbb R^N$ of radius $R$ centered at the origin
and  let $B$ and $\Gamma$ denote the Beta and the Gamma function, respectively.
Our main  result, contained in Section 5,  is the following.
\begin{theorem}
\label{maintheorem}
Let $N\in \mathbb{N} $, $N\geq 2$, 
$k,l \in \mathbb{R} $, $\alpha > 0$ and $l+N+\alpha >0$. 
Further, assume that one of the following conditions holds:
\\
{\bf (i)} $l+1\leq k $;
\\
{\bf (ii)} $k\leq l+1$ and $ l\frac{N+\alpha-1}{N+\alpha} \leq k\leq 0$; 
\\  
{\bf (iii)} $N\geq 2$, $ 0\leq k\leq l+1$ and 
\begin{equation}\label{l_1N3}
l\le l_1 (k,N,\alpha )  :=  \frac{(k+N+\alpha-1)^3 }{(k+N+\alpha-1)^2 - \frac{(N+\alpha-1)^2 }{N+\alpha} } -N -\alpha\,.
\end{equation}
\\
Then  
\begin{equation}
\label{mainineq}
\dint_{\partial \Omega } |x|^k x_N^\alpha\, {\mathcal H}_{N-1} (dx)
\geq 
C_{k,l,N, \alpha} ^{rad}  
\left( 
\dint_{\Omega } |x|^lx_N^\alpha\, dx 
\right) 
^{(k+N+\alpha-1)/(l+N+\alpha) } , 
\end{equation}
for all smooth sets $\Omega $ in 
$\mathbb{R}^N _+  $,
where 

\begin{eqnarray}
\label{defCkl}
C_{k,l,N, \alpha} ^{rad} & := &
\frac{\dint_{\partial B_1 } |x|^k x_N^\alpha\, {\mathcal H}_{N-1} (dx)}
{\left( \dint_{B_1 \cap  \mathbb{R}^N _+}
 |x|^l x_N^\alpha\, dx \right ) ^{(k+N+\alpha-1)/(l+N+\alpha) } } 
\\
 &= &
\nonumber
\left( l+\alpha +N\right) ^{\frac{k+N+\alpha -1}{l+N+\alpha }}\left(
B\left( \frac{N-1}{2},\frac{\alpha +1}{2}\right) \frac{\pi ^{\frac{N-1}{2}}}{
\Gamma \left( \frac{N-1}{2}\right) }\right) ^{\frac{l-k+1}{l+N+\alpha }}.
\end{eqnarray}
Equality in (\ref{mainineq}) holds  if $\Omega =B_R\cap\mathbb{R}_+^N$.
\end{theorem}

\noindent

Note that the weights we consider are not radial and it seems not trivial to use
 spherical symmetrization.  So that we did not try to   adapt the techniques contained in  \cite{diGiosia_etal},
 and, depending on the regions where the three parameters lie, we use different methods.
The proof in the case {\bf (i)} is given in \cite{ABCMP_atti}. It is based on Gauss's Divergence Theorem.
In the case {\bf (ii)}  (see Theorem \ref{th1bis}) the proof uses an appropriate  change of variables, 
which has been introduced in \cite{H} and \cite{HK},
together with the isoperimetric inequality with respect to the weight $x_{N}^{\alpha} $.
 The case {\bf (iii)}  (see Theorem \ref{th1ter}) is the most delicate and it requires several different arguments:
 again a suitable change of variables,  then an interpolation argument, introduced for the first time in our previous paper \cite{ABCMP} and, finally, the so-called starshaped rearrangement.

In Section 4 we provide some necessary conditions on $k$, $l$ and $\alpha$  such that 
the half-ball centered at the origin is an isoperimetric set. In the proof we 
 firstly evaluate the second variation of the perimeter functional.
 The claim is achieved using the  fact that such  a variation 
 at a minimizing set  must be nonnegative, together with a nontrivial 
weighted Poincar\' e inequality on the sphere derived in \cite{BCM2}.

\noindent Part of these results have been announced in \cite{ABCMP_atti}.

\section{Notation and preliminary results  } 

Throughout this article  $N$ will denote a natural number with $N\geq 2$,  
$k$ and $l $ are real numbers, while $\alpha$ is a nonnegative number
and 
\begin{equation}
\label{ass1}
 l+N+\alpha>0 .
\end{equation}

\noindent 
Let us introduce some notation. 
\begin{eqnarray*}
\mathbb{R}^{N}_{+}
& := & 
\left\{ x \in \mathbb{R}^{N}: x_N  >0 \right\}, 
\\
\mathbb{S}^{N-1}_{+}
& := & 
\left\{ x \in \mathbb{S}^{N-1} : x_N  >0 \right\},
\\
B_{R}(x_0 )
 & := &
\left\{ x\in \mathbb{R}^{N}:\left\vert x-x_0 \right\vert <R\right\} , \quad (x_0 \in \mathbb{R}^N ),
\\
B_{R} 
& := & 
B_{R}(0), \quad (R>0),
\\
B_{R}^{+}  
& := & 
B_{R}  \cap \mathbb{R}^{N}_{+}.
\\
\end{eqnarray*}
Furthermore, ${\mathcal L} ^m $ will denote  the $m$-dimensional Lebesgue measure, ($1\leq m\leq N$), and
\begin{eqnarray*} 
\omega _N & := & {\mathcal L}^N (B_1 ),
\\
\kappa (N, \alpha ) & := & {\mathcal L} ^{N-1} (\mathbb{S} ^{N-1} _+ ).
\end{eqnarray*}
Note that 
\begin{equation}
\label{measSN-1+}
\kappa (N, \alpha ) = B\left( \frac{
N-1}{2},\frac{\alpha +1}{2}\right) 
\frac{\pi ^{\frac{N-1}{2}}}{\Gamma \left( 
\frac{N-1}{2}\right) },
\end{equation}
where $B$ and $\Gamma$ are the Beta function and the Gamma function,
respectively,
(see \cite{BCM3}).
\\
We will use frequently $N$-dimensional spherical coordinates $(r, \theta)$ in $\mathbb{R} ^N$: 
$$
\mathbb{R}^N \ni x = r\theta , \quad \mbox{where $r=|x|$, and $\theta = x|x|^{-1} \in \mathbb{S}^{N-1} $.}
$$ 
If $M$ is any set in $\mathbb{R}^N _{+}$, then $\chi _M $ will denote its characteristic function.

\noindent Next, let $k$ and $l$ be real
numbers satisfying (\ref{ass1}). We define a measure $\mu _{l, \alpha}$ by 
\begin{equation}
d\mu _{l, \alpha}(x)=|x|^{l} x_N^\alpha\,dx.  
\label{dmu}
\end{equation}
If $M \subset $ ${\mathbb R}^{N}_{+}$ is a  measurable set with finite 
$\mu _{l, \alpha} $-measure, then we define $M^{\star }$, the 
\\
$\mu_{l,\alpha }$-symmetrization of $M$,
 as follows:
\begin{equation}
M^{\star } := B_{R}^{+} 
\hspace{.2 cm}
 \text{with }
 R:
\mu_{l, \alpha} 
\left( B_{R}^{+}    \right) =
\mu _{l, \alpha} \left( M\right) = \int_M d\mu _{l, \alpha} (x) .  
\label{mu_(M)}
\end{equation}
If $u: \mathbb{R}^{N}_{+} \rightarrow \mathbb{R} $ is a measurable function such that
$$
\mu_{l, \alpha} \left( \left\{ |u(x)|>t\right\} \right) <\infty \qquad \forall t>0,
$$
then let $u^{ \star }$ denote the weighted Schwarz symmetrization of $u$, or,
in short, the \\
$\mu_{l, \alpha} -$symmetrization of $u$, which is given by
\begin{equation}
u^{ \star }(x)=\sup \left\{ t\geq 0:\mu_{l, \alpha}
 \left( \left\{ |u(x)| >t\right\} \right)
>
\mu _{l, \alpha} \left( B_{\left\vert x\right\vert }^{+}   \right) \right\} .  
\label{u_star}
\end{equation}

Note that $u^{\star }$ is radial and radially non-increasing, 
and if $M$ is a measurable set with finite $\mu _l $-measure, then 
$$
\left( \chi _M \right) ^{\star} = \chi _{M^{ \star }} .
$$

The {\sl $\mu _{k, \alpha}$--perimeter\/} of a measurable set $M $ is given by 
\begin{equation}
P_{\mu _{k, \alpha}}(M ):=\sup \left\{ \int_{M }\mbox{div}\,\left(x_N^\alpha |x|^{k}
\mathbf{v}\right) \,dx:\,\mathbf{v}\in C_{0}^{1}(\mathbb{R}^N ,\mathbb{R}^{N}),\,|
\mathbf{v}|\leq 1\mbox{ in }\, M \right\} .
\end{equation}

\noindent It is well-known that the above 
\textsl{distributional definition} of
weighted perimeter is equivalent to the following

\begin{equation}
P_{\mu_{k} }(M )
= \left\{ 
\begin{array}{ccc}
\displaystyle\int_{\partial \Omega }|x|^{k} \, {\mathcal H} _{N-1}(dx)
 & \mbox{ if } & 
\partial \Omega  \mbox{ is } (N-1)-\mbox{rectifiable } \\ 
&  &  \\ 
+ \infty \qquad
 & \mbox{ otherwise,} & 
\end{array}
\right.
\end{equation}
where, here and throughout, ${\mathcal H} _{N-1} $ will denote the  $(N-1)$-dimensional Hausdorff-measure. 

We will call a set $\Omega \subset \mathbb{R}^N _+ $ {\sl smooth}, 
if 
 for every $x_0 \in \partial \Omega \cap \mathbb{R}^N _+ $, 
there is a number $r >0$ such that $B_r (x_0 )\subset \mathbb{R}^N _+ $,  $B_r (x_0 ) \cap \Omega $
has exactly one connected component and $B_r (x_0 ) \cap \partial \Omega $ 
is the graph of a $C^1 $--function on an open set in $\mathbb{R} ^{N-1} $.  

Let $\Omega \subset \mathbb{R} ^{N}_{+}$ and $p\in \left[ 1,+\infty \right) $.
We will denote by $L^{p}(\Omega ,d\mu
_{l, \alpha})$ the space of all Lebesgue measurable real valued functions $u$ such
that
\begin{equation}
\left\Vert u \right\Vert 
_{L^{p}(\Omega ,d\mu _{l, \alpha})}
:=\left( 
\int_{\Omega
}
\left\vert 
u\right\vert 
^{p}d\mu _{l, \alpha} (x) 
\right) ^{1/p}
<+\infty .
\label{Norm_Lp}
\end{equation}
\\  
By $W^{1,p}(\Omega ,d\mu _{l, \alpha})$ we denote the weighted Sobolev space
consisting of all functions which together with their weak derivatives $u_{x_{i}}$, ($i=1,...,N$), 
belong to $L^{p}(\Omega ,d\mu _{l, \alpha})$.
This space will be equipped with the norm
\begin{equation}
\left\Vert u\right\Vert _{W^{1,p}(\Omega ,d\mu _{l, \alpha})}:=\left\Vert
u\right\Vert _{L^{p}(\Omega ,d\mu _{l, \alpha})}+\left\Vert \nabla u\right\Vert
_{L^{p}(\Omega ,d\mu _{l, \alpha})}.  
\label{Norm_Wp}
\end{equation}
Finally,  ${\mathcal D} ^{1,p}( \Omega  ,d\mu _{k, \alpha})$ will stand for the closure of 
$C_{0}^{\infty }(\mathbb{R}^N )$ under the norm 
$$
\left( \int_{\Omega } |\nabla u|^p \, d\mu _{k,\alpha } (x) \right) ^{1/p}.
$$
We will often use the following well-known {\sl Hardy-Littlewood inequality} 
\begin{equation}
\label{hardylitt1}
\int_{\mathbb{R}^{N}_{+} } uv \, d\mu _{l, \alpha}(x) \leq \int_{\mathbb{R}^{N}_{+} } u^{ \star} v^{\star} \, d\mu _{l, \alpha} (x) ,
\end{equation}
which holds for any couple of functions 
$u,v\in L^2 (\mathbb{R}^{N}_{+} , d\mu _{l, \alpha} )$.

Now let us recall the so-called starshaped rearrangement 
(see \cite{Kaw}) which we will use in Section 5.
For later convenience, we will write $y$ for points in $\mathbb{R}^{N}_{+} $  and  $(z, \theta )$ for corresponding $N$-dimensional spherical coordinates ($z= |y|$, $\theta = y|y|^{-1} $).
\\
We call a measurable set $M\subset \mathbb{R}^{N}_{+} $ {\sl starshaped\/} if the set 
$$
M\cap \{ z\theta : \, z\geq 0 \}
$$ 
is either empty or a segment 
$\{ z\theta : \, 0\leq z< m(\theta ) \} $ 
for some number $m(\theta ) >0 $, for almost every $\theta \in {\mathbb S} ^{N-1} $. 
\\
If $M$ is a bounded measurable set in $\mathbb{R}^{N}_{+} $, and 
$\theta \in {\mathbb S}^{N-1}_{+} ,$ then let
$$
M(\theta ) := M\cap \{ z\theta :\, z\geq 0\}.  
$$
There is a unique number $m(\theta )\in [0,+\infty )$ such that
$$
\int_0 ^{m(\theta )} z^{N-1}\, dz = \int_{M(\theta )} z^{N-1} \, dz.
$$
We define 
$$
\widetilde{M}(\theta ) := \{ z\theta : \, 0\leq z\leq m(\theta )  \} ,
\quad (\theta \in {\mathbb S} ^{N-1}_{+} ),
$$
and 
$$
\widetilde{M} := \{ z\theta : \, z\in \widetilde{M}(\theta ) , \, 
\theta \in {\mathbb S} ^{N-1}_{+} \} .
$$
We call the set $\widetilde{M}$ the {\sl starshaped rearrangement of $M$\/}. 
\\
Note that $\widetilde{M} $ is Lebesgue measurable and starshaped, and 
we have
\begin{equation}
 \label{starsh1}
 {\mathcal L} ^N (M) = {\mathcal L} ^N (\widetilde{M}).
\end{equation}
If $v:\mathbb{R}^{N}_{+} \to \mathbb{R} $ is a measurable function with compact support, 
and $t\geq 0$, then let
$ E_t $ be the super-level set $\{ y: \, |v(y)| \geq t\} $. We
define 
$$
\widetilde{v} (y) := \sup \{ t\geq 0 :\, y \in \widetilde{E_t }   \} .
$$
We call $\widetilde{v} $ the {\sl starshaped rearrangement of $v$ \/}.  
It is easy to verify that 
$\widetilde{v}$ is equimeasurable with $v$, 
that is, the following properties hold:
\begin{eqnarray}
\label{starsh2}
 & & \widetilde{E_t} = \{ y:\, \widetilde{v} (y)\geq t\} , 
 \\
 \label{starsh3}
 & & {\mathcal L} ^N (E_t ) = 
 {\mathcal L} ^N (\widetilde{E_t} ) \quad \forall t\geq 0.
\end{eqnarray}
This also implies  Cavalieri's principle: 
If $F\in C ([0, +\infty ))$ with $F(0)=0$ and if 
$F(v) \in L^1 ( \mathbb{R}^N ) $, then
\begin{equation}
\label{caval1}
\int_{\mathbb{R} ^N } F(v)\, dy =  
\int_{\mathbb{R} ^N } F(\widetilde{v} )\, dy
\end{equation}
and if $F$ is non-decreasing, then 
\begin{equation}
\label{monrearr}
\widetilde{F(v)} = F(\widetilde{v}).
\end{equation}
Note that the mapping
$$
z\longmapsto \widetilde{v} (z\theta ) , \quad (z\geq 0),
$$
is non-increasing for all $\theta \in {\mathbb S} ^{N-1} $.
\\ 
If $v,w\in L^2 (\mathbb{R}^{N}_{+} )$ are functions with compact support,
then there holds 
Hardy-Littlewood's inequality:
\begin{equation}
\label{harlit}
\int_{\mathbb{R}^{N}_{+} } vw \, dy \leq  \int_{\mathbb{R}^{N}_{+} } 
\widetilde{v} \widetilde{w} \, dy.
\end{equation}    
If $f:(0,+\infty) \to \mathbb{R} $ is a measurable function with compact support, then its (equimeasurable) {\sl non-increasing rearrangement }, $\widehat{f} : (0,+\infty )\to [0,+\infty )$, 
is the monotone non-increasing 
function such that 
$$
{\mathcal L} ^1 \{ t  \in [0,+\infty )  :\,  |f(t )| > c\}  = 
{\mathcal L}^1 \{ t \in [0,+\infty ) :\,  \widehat{f}(t ) > c \}  \quad \forall c\geq 0,
$$
see \cite{Kaw}, Chapter 2.
A general P\'{o}lya-Szeg\"o principle for non-increasing rearrangement has been given in \cite{Lan}, Theorem 2.1. For later reference we will only need a special case: 
\begin{lemma}
\label{Landes} 
Let $\delta \geq 0$,  and let $f:(0,+\infty ) \to \mathbb{R} $ be a bounded, locally Lipschitz continuous function with bounded support, such that 
$$
\int_0 ^{+\infty } t  ^{\delta } |f' (t ) |\, dt <+\infty .
$$
Then $\widehat{f} $ is locally Lipschitz continuous and 
\begin{equation}
\label{landes1}
\int_0 ^{+\infty } t^{\delta } |\widehat{f}' (t ) |\, dt
\leq
\int_0 ^{+\infty } t^{\delta } |f' (t ) |\, d t.
\end{equation}  
\end{lemma}
 \bigskip

\section{The functionals ${\mathcal R}_{k,l,N,\alpha}$ and ${\mathcal Q}_{k,l,N,\alpha}$ }
Throughout this section  we  assume (\ref{ass1}), i.e.
\begin{equation*}
k+N+\alpha-1 >0 \ \mbox{ and } \ l+N+\alpha>0 .
\end{equation*}

If $M $ is any measurable subset of $\mathbb R^{N}_{+}$, with $0<\mu _{l,\alpha} (M)<+\infty $, we set
\begin{equation}
\label{rayl1}
{\mathcal R}_{k,l,N, \alpha} (M) := 
\frac {P_{ \mu_{k , \alpha}}  (M) }
{  \left( \mu_{l,\alpha} (M) \right)^{(k+N+\alpha-1)/(l+N+\alpha)} }. 
\end{equation}
Note that
\begin{equation}
\label{Rklsmooth}
{\mathcal R} _{k,l,N,\alpha} (M ) = 
\frac{
\dint_{\partial M }x_N^\alpha |x|^k \, {\mathcal H}_{N-1} (dx) 
}{
\left( \dint_{M }x_N^\alpha |x|^l \, dx \right) ^{(k+N+\alpha-1)/(l+N+\alpha)} } 
\end{equation}
if the set $M$ is smooth.

If $u\in C_0 ^1  (\mathbb{R} ^N_+ )\setminus \{ 0\} $, we set
\begin{equation}
\label{rayl2}
{\mathcal Q}_{k,l,N, \alpha} (u ) := \frac{\dint_{\mathbb{R} ^N_+ }x_N^\alpha |x|^k |\nabla u| \, dx}{ 
\left( \dint_{\mathbb{R} ^N_+ } x_N^\alpha|x|^l |u| ^{(l+N+\alpha)/(k+N+\alpha-1)} \, dx \right) ^{(k+N+\alpha-1)/(l+N+\alpha)}}. 
\end{equation}

Finally, we define 
\begin{equation}
\label{isopco}
C_{k,l,N, \alpha}^{rad} := {\mathcal R}_{k,l,N, \alpha}(B_1 \cap {\mathbb{R} ^N_+ }).
\end{equation} 
We study the following isoperimetric problem: 
\\[0.5cm]
{\sl Find the constant $C_{k,l,N, \alpha} \in [0, + \infty )$, such that}
\begin{equation}
\label{isopproblem}
C _{k,l,N, \alpha} := \inf \{ {\mathcal R}_{k,l,N, \alpha} (M):\, 
\mbox{{\sl $M$ is measurable with $0<\mu _{l,\alpha} (M) <+\infty $.}} \} 
\end{equation}   
Moreover, we are interested in conditions on $k$, $l$ and $\alpha$ such that
\begin{equation}
\label{isoradial}
{\mathcal R}_{k,l,N, \alpha} (M) \geq {\mathcal R}_{k,l,N, \alpha} (M^{ \star} ) 
\end{equation}
holds for all measurable sets $M\subset {\mathbb{R} ^N_+ }$ with  $ 0<\mu _{l,\alpha}(M)<+\infty $. 
\\[0.1cm]
Let us begin with some immediate observations.
\\
If $M$ is a measurable subset of $\mathbb R^{N}_{+}$ with finite $\mu _{l,\alpha}$-measure and $\mu_{k,\alpha} $-perimeter,
 then there exists a sequence of smooth sets 
$\{ M_n \} $
 such that 
$$\lim_{n\to \infty } \mu _{l,\alpha} (M_n \Delta M) =0   \,\,\,\,  \text{and} \,\,
\lim_{n\to \infty } P_{\mu _{k,\alpha} } (M_n ) = P_{\mu _{k,\alpha}} (M) .
$$
This property is well-known for Lebesgue measure (see for instance 
\cite{G}, Theorem 1.24) 
and its proof carries over to the weighted case. This implies that we also have 
\begin{equation}
\label{CklNsmooth}
C_{k,l,N, \alpha} = \inf \{ {\mathcal R}_{k,l,N, \alpha} (\Omega ):\, \Omega \subset \mathbb{R} ^{N}_+, \, \Omega  
\mbox{ smooth} \} .
\end{equation}  
The functionals ${\mathcal R}_{k,l,N, \alpha } $ and ${\mathcal Q}_{k,l,N,\alpha } $ 
have the following homogeneity properties, 
\begin{eqnarray}
\label{hom1}
 {\mathcal R}_{k,l,N, \alpha } (M ) & = & {\mathcal R}_{k,l,N, \alpha } (tM ) ,
\\
{\mathcal Q}_{k,l,N, \alpha } (u) & = & {\mathcal Q}_{k,l,N,\alpha } (u^t ),
\end{eqnarray}
where  $t>0$, $M $ is a measurable set with $0<\mu_{l, \alpha} (M)<+\infty $, 
$u\in C_0 ^1 (\mathbb{R}^N_+ )\setminus \{ 0\}$,  
\\
$tM := \{tx:\, x\in M \} $ 
and $u^t (x):= u(tx) $, ($x\in \mathbb{R} ^N_+ $), and there holds   
\begin{equation}
\label{isopconst2}
C_{k,l,N, \alpha} ^{rad} = {\mathcal R}_{k,l,N, \alpha} (B_1^{+} ).
\end{equation}  
Hence we have that 
\begin{equation}
\label{relCC}
C_{k,l,N, \alpha} \leq C_{k,l,N, \alpha} ^{rad} ,
\end{equation}  
and (\ref{isoradial}) holds if and only if 
$$
C_{k,l,N,\alpha } = C_{k,l,N,\alpha } ^{rad} .
$$  
Finally, we recall the following weighted isoperimetric inequality proved, for example, in \cite{BCM2}  (see also \cite{XR} and \cite{ MadernaSalsa}).

\begin{proposition}\label{BCM2} 
For all measurable sets $M\subset \mathbb{R} ^N_+$, with $0< \mu _{0, \alpha}  (M)<+\infty $, the following inequality holds true
\begin{equation}\label{isopclass}
{\mathcal R} _{0,0,N, \alpha} (M) := 
\frac {P_{ \mu_{0, \alpha}}  (M) }
{  \left( \mu_{0,\alpha} (M) \right)^{(N+\alpha-1)/(N+\alpha)} } \geq C_{0,0,N, \alpha} ^{rad}:= 
\frac {P_{ \mu_{0, \alpha}}  (M^{ \star}) }
{  \left( \mu_{0,\alpha} (M^{ \star}) \right)^{(N+\alpha-1)/(N+\alpha)} } \,,
\end{equation}
where $M^{\star}=B_{R}^{+} $ with $R$ such that $\mu_{0, \alpha}(M)=\mu_{0, \alpha}(M^{ \star})$
\end{proposition}

We recall that the isoperimetric constant $C_{0,0,N, \alpha} ^{rad}$ is explicitly computed in \cite{BCM2}, see also \cite{MadernaSalsa} for the case $N=2$.

\begin{lemma}
\label{hardylitt}
Let $l>l' >-N -\alpha$. Then 
\begin{equation} 
\label{hardylitt2}
\frac{\left( \mu _{l, \alpha} (M)  \right) 
^{1/(l+N+\alpha)} 
}{ 
\left( \mu _{l', \alpha} (M)  \right)
^{1/(l'+N+\alpha)}
} 
\geq \frac{\left( \mu _{l, \alpha} (M^{ \star})  \right) 
^{1/(l+N+\alpha)} 
}{ 
\left( \mu _{l', \alpha} (M^{ \star})  \right)
^{1/(l'+N+\alpha)}
} 
\end{equation}
for all measurable sets $M\subset \mathbb{R} ^N_+ $ with  $0<\mu_{l, \alpha}(M)<+\infty $.
Equality holds only for half-balls $B_{R}^{+} $, ($R>0$).
\end{lemma}

{\sl Proof: } Let $M^{ \star} $ be the 
$\mu_{l, \alpha} $-symmetrization of $M$. Then we obtain, using the Hardy-Littlewood inequality,
\begin{eqnarray*}
\mu _{l' , \alpha} (M) =\int_Mx_N^\alpha |x| ^{l'} \, dx & = & \int_{\mathbb{R} ^N_+ }  |x|^{l'-l} \chi _M (x)\, d\mu _{l, \alpha} (x) 
\\
 & \leq & 
\int_{\mathbb{R} ^N _+}  \left( |x|^{l'-l} \right) ^{ \star}  \left( \chi _M  \right) ^{ \star} (x)\, d\mu _{l, \alpha} (x) 
\\
 & = & 
\int_{\mathbb{R} ^N_+ }   |x|^{l'-l}  \chi _{M^{ \star} }  (x)\, d\mu _{l, \alpha} (x) 
\\
 & = & 
\int_{M^{ \star} } x_N^\alpha |x|^{l' }\, dx =\mu _{l', \alpha } (M^{ \star} ).
\end{eqnarray*} 
This implies (\ref{hardylitt2}). 

\noindent Next assume that equality holds in (\ref{hardylitt2}). Then we must have 
$$
\int_M |x|^{l'-l} \, d\mu _{l, \alpha} (x) = \int_{M^{ \star} } |x|^{l'-l} d\mu _{l, \alpha} (x) ,
$$
that is, 
$$
\int_{M\setminus M^{ \star} } |x|^{l'-l} \, d\mu _{l, \alpha} (x) = \int_{M^{ \star} \setminus M} |x|^{l'-l} d\mu _{l, \alpha}  (x) .
$$
Since $l'-l<0$, this means that 
$ \mu _l ( M\Delta M^{ \star} )=0$. The Lemma is proved. 
$\hfill \Box $
\begin{lemma}
\label{rangekl1}
Let $k,l, \alpha$ satisfy (\ref{ass1}). Assume that $l>l' >-N-\alpha$ and 
$C_{k,l,N, \alpha}  = C_{k,l,N, \alpha} ^{rad} $. 
Then we also have 
$C_{k,l',N, \alpha}  = C_{k,l',N, \alpha} ^{rad} $.
Moreover, if  $
{\mathcal R}_{k,l',N, \alpha} (M ) = C_{k,l',N, \alpha} ^{rad} $  
for some measurable set $M\subset \mathbb{R} ^N_+ $, with $0< \mu _{l' , \alpha} (M) <+\infty $,
then $M = B_{R}^{+}$ for some $R>0$.
\end{lemma}
{\sl Proof:} By our assumptions and Lemma \ref{hardylitt}  we have for every measurable set $M$ with 
$0<\mu _{l, \alpha}(M) <+\infty $,  
\begin{eqnarray*}
{\mathcal R}_{k,l',N, \alpha} (M ) & = & {\mathcal R}_{k,l,N, \alpha} (M ) 
\cdot 
\left[
\frac{
\left( 
\mu _{l, \alpha} (M) 
\right) ^{1/(l+N+\alpha)}
}{ 
\left( \mu_{l', \alpha} (M) 
\right) ^{1/(l'+N+\alpha)}
} 
\right] ^{k+N+\alpha-1}
\\
 & \geq & 
C_{k,l',N, \alpha}^{rad},
\end{eqnarray*}
with equality only if 
$M = B^{+}_{R}$ for some $R>0$. 
$\hfill \Box $

\begin{lemma}
\label{R2} 

Assume that $k \leq l+1$. Then
\begin{equation}
\label{ineqQR}
C_{k,l,N, \alpha} 
=
\inf \left\{ {\mathcal Q}_{k,l,N, \alpha} (u) :\, u\in C_0 ^1 (\mathbb{R}_+ ^N 
 )\setminus \{ 0\} \right\} . 
\end{equation}
\end{lemma}

{\sl Proof: }
The proof uses classical arguments (see, e.g. \cite{FleRi}). 
We may restrict ourselves to nonnegative functions $u$.  
By (\ref{isopproblem}) and the coarea formula we obtain,
\begin{eqnarray}
\label{coarea1}
\int_{\mathbb{R} ^N_+ }x_N^\alpha |x|^k |\nabla u| \, dx & = & 
\int _0 ^{\infty } \int\limits_{u=t } x_N^\alpha |x|^k \, {\mathcal H} _{N-1 } (dx) \, dt 
\\
\nonumber
 & \geq & C_{k,l,N, \alpha} \int_0 ^{\infty } \left( \int_{u>t } x_N^\alpha |x|^l \, dx 
\right) ^{(k+N+\alpha-1)/(l+N+\alpha)} \, dt.
\end{eqnarray}
Further,
Cavalieri's principle gives
\begin{equation}
\label{cavalieri}
u(x)= \int_0 ^{\infty } \chi _{\{ u>t\} } (x)\, dt , \quad (x\in \mathbb{R} ^N ).
\end{equation}
Hence (\ref{cavalieri}) and Minkowski's inequality for integrals (see \cite{Stein}) lead to 
\begin{eqnarray}
\label{ineqmeas}   
 & & 
\\
\nonumber  
&&\int_{\mathbb{R} ^N_+ }x_N^\alpha |x|^l |u|^{(l+N+\alpha)/(k+N+\alpha-1)}  \, dx \qquad \qquad  \\
\nonumber  &=&  \int_{\mathbb{R} ^N _+}x_N^\alpha |x|^l \left| \int_0 ^{\infty} 
\chi_{\{ u>t\} } (x)\, dt \right| ^{(l+N+\alpha)/(k+N+\alpha-1)} \, dx 
\\
\nonumber 
 & \leq & \left( \int_0 ^{\infty } \left( 
\int_{\mathbb{R}^N_+ }x_N^\alpha |x|^l \chi _{\{ u>t \} } (x) \, dx \right) 
^{(k+N+\alpha-1)/(l+N+\alpha)} \, dt \right) ^{(l+N+\alpha)/(k+N+\alpha-1)} 
\\
\nonumber
 & = & \left( \int_0 ^{\infty } \left( \int_{ u>t }x_N^\alpha |x|^l \, dx \right) 
^{(k+N+\alpha-1)/(l+N+\alpha)} dt \right) ^{(l+N+\alpha)/(k+N+\alpha-1)} .
\end{eqnarray}
Now (\ref{coarea1}) and (\ref{ineqmeas}) yield
\begin{equation}
\label{ineqQ1}
{\mathcal Q}_{k,l,N, \alpha} (u) \geq C_{k,l,N, \alpha} \quad \forall u\in C_0 ^1\setminus \{ 0\}  (\mathbb{R}_+ ^N ).
\end{equation}
To show (\ref{ineqQR}), 
let $\varepsilon >0$, and choose a smooth set 
$\Omega $ such that 
\begin{equation}
\label{ineqR1}
{\mathcal R}_{k,l,N,\alpha} (\Omega ) \leq C_{k,l,N,\alpha } +\varepsilon .
\end{equation}
It is well-known that there exists a sequence $\{ u_n \} \subset 
C_0 ^{\infty } (\mathbb{R} ^N )\setminus \{ 0\} $ 
such that 
\begin{eqnarray}
\label{limperim}
\lim_{n\to \infty } \int_{\mathbb{R}^N _+} x_N^\alpha |x|^k |\nabla u_n | \, dx = 
\int_{\partial \Omega }x_N^\alpha |x|^k \, {\mathcal H} _{N-1} (dx) ,
\\
\label{limmeas}
\lim_{n\to \infty } \int_{\mathbb{R}_+ ^N } x_N^\alpha |x|^l 
|u_n |^{(l+N+\alpha)/(k+N+\alpha-1)}  \, dx = \int_{ \Omega } x_N^\alpha|x|^l \, dx.
\end{eqnarray}
To do this, one may choose mollifiers of $\chi _{\Omega } $ 
as $u_n $ (see e.g. \cite{Talenti1}).  
Hence, for large enough $n$ we have
\begin{equation}
\label{ineqQ2}
{\mathcal Q}_{k,l,N,\alpha } (u_n ) \leq C_{k,l,N,\alpha } + 2\varepsilon .
\end{equation}
Since $\varepsilon $ was arbitrary, (\ref{ineqQR}) now 
follows from (\ref{ineqQ1}) and (\ref{ineqQ2}).      
$\hfill \Box $

\section{Necessary conditions}

In this section we assume that 
\begin{equation*}
k+N+\alpha-1 >0 \ \mbox{ and } \ l+N+\alpha>0 .
\end{equation*}
The main result is Theorem \ref{R4}  which highlights
the phenomenon of symmetry breaking.

\noindent The following result holds true.

\begin{lemma}
\label{R3}
A necessary condition for 
\begin{equation}
\label{C>0}
C_{k,l,N, \alpha} >0
\end{equation}
 is 
\begin{equation}
\label{k_l_ineq1}
 l \frac{N+\alpha-1}{N+\alpha} \leq k .
\end{equation}
\end{lemma}

{\sl Proof:} Assume that $k<l(N+\alpha-1 )/(N+\alpha)$, and let 
$te_1 = (t, 0, \ldots , 0)$, ($t>2$). Since for any $x\in B_1 (te_1 )$,  it results $t-1\le |x|\le t+1$, we have 
$$
{\mathcal R}_{k,l,N, \alpha} (B_1 (te_1 ) ) \leq 
D \frac{ (t +1)^k}{ (t-1) ^{l (k+N+\alpha-1)/(l+N+\alpha)} }.
$$
where the positive constant $D= D(k,l, N, \alpha) $ is given by
$$
D=\frac{
\dint_{\partial (B_1 (te_1 )\cap \mathbb{R}_+^{N})}x_N^\alpha \, {\mathcal H}_{N-1} (dx) 
}{
\left( \dint_{B_1 (te_1 )\cap \mathbb{R}_+^{N} }x_N^\alpha \, dx \right) ^{(k+N+\alpha-1)/(l+N+\alpha)} } 
$$
Since $k-l(k+N+\alpha-1)/(l+N+\alpha) <0$, it follows that
$$
\lim_{t\to \infty } {\mathcal R}_{k,l,N,\alpha } (B_1 (te_1 ) ) =0.
$$
$\hfill \Box $

\begin{theorem}
\label{R4} 

A necessary condition for 
\begin{equation}
\label{isop1}
 C_{k,l,N, \alpha} = C_{k,l,N, \alpha} ^{rad} 
\end{equation}
is 
\begin{equation}
\label{k_l_ineq2}
l+1 \leq k + \frac{N+\alpha-1}{ k+N+\alpha-1} .
\end{equation}
\end{theorem}
\medskip

\begin{remark} \rm
Theorem \ref{R4} means that if  $l+1 \leq k + \frac{N+\alpha-1}{ k+N+\alpha-1}$, 
 then symmetry breaking occurs, that is $C_{k,l,N, \alpha} < C_{k,l,N, \alpha} ^{rad} $. 
Our proof relies on the fact that the second variation of the perimeter for smooth volume-preserving 
perturbations from the ball $B_{1}^{+} $ is non-negative if and only if (\ref{k_l_ineq2}) holds. Note that this also 
follows from a general second variation formula with  volume and perimeter densities, see \cite{Mo2}. 
\end{remark}

{\sl Proof:} First we assume $N\geq 2$. Let $(r, \theta )$ denote 
$N$--dimensional spherical coordinates, such that
$$
\theta _1 = \arccos \frac{x_N}{|x|} , \quad\theta_1 \in [0, \pi /2],
$$
and 
$u \in C^2 (\mathbb{S}_{+}^{N-1} )$, $s\in C^2 (\mathbb{R})$ with $s(0)=0$, 
and define
$$
U(t ) := \{ x=r\theta \in \mathbb{R}_+ ^N : \, 0\leq r < 1+ t u(\theta ) + s(t) \} , 
\quad (t\in \mathbb{R} ).
$$
Note that $U(0)= B_{1}^{+} $. 
By the Implicit Function Theorem, we may choose $s$ in such a way that 
\begin{equation}
\label{intid1}
\int_{U(t)}x_N^\alpha |x|^l \, dx = \int _{B_{1}^{+} } x_N^\alpha |x|^l \, dx \quad \mbox{for $|t|<t_0$},
\end{equation}
for some number $t_0 >0$. We set $s_1 := s'(0) $ and $s_2 := s^{\prime \prime} (0)$. 
Let $d \Theta$ be the surface element on the sphere
and 
\begin{equation}
\label{h}
h:= h(\theta_1)  = \cos^{\alpha} \theta_1 = \left(  \frac{x_N}{|x|} \right)^{\alpha}.
\end{equation}
Since
$$
\int_{U(t)} x_N^\alpha|x|^l \, dx = \int_{\mathbb{S}_+ ^{N-1 } } h
\int_0 ^{1+ t u(\theta ) + s(t)}  \rho ^{l+N+\alpha-1} \, d\rho \, d\Theta,
$$
a differentiation at $t=0$ of (\ref{intid1}) leads to
\begin{eqnarray}
\label{intid2}
0 & = & \int_{\mathbb{S}_+ ^{N-1} }  (u+ s_1 )\, h d\Theta \quad \mbox{and }
\\
\label{intid3}
0 & = & (l+N+\alpha-1) \int_{\mathbb{S}_+^{N-1} }  (u+ s_1 )^2  h \, d\Theta + s_2 
\int_{\mathbb{S}_+^{N-1} } h \,  d\Theta .
\end{eqnarray}
Next we consider the perimeter functional 
\begin{eqnarray}
\label{perim}
J(t) & := & \int_{\partial U(t)} x_N^\alpha |x|^k \, {\mathcal H}_{N-1} (dx)  
\\
\nonumber
 & = & \int_{\mathbb{S}_+ ^{N-1} } (1+tu + s(t) )^{k+N+\alpha-2} 
\sqrt{ (1+ tu+s(t) )^2 + t^2 |\nabla _{\theta } u|^2 } \,  h \, d\Theta ,
\end{eqnarray}
where $\nabla _{\theta }$ denotes the gradient on the sphere. 
Differentiation at $t=0$ of (\ref{perim}) leads to 
\begin{eqnarray*}
J'(0) & = & (k+N+\alpha-1) \int_{\mathbb{S} _+^{N-1} } (u+ s_1 ) \, h \, d\Theta , \quad \mbox{and }
\\
J^{\prime \prime } (0) & = & 
(k+N+\alpha-2) (k+N+\alpha-1) \int_{\mathbb{S} _+^{N-1} }  (u+s_1 )^2  \, h \, d\Theta + 
\\
 & & + (k+N+\alpha-1) s_2 \int_{\mathbb{S} _+^{N-1} } \, h \, d\Theta +  \int_{\mathbb{S}_+ ^{N-1} } 
 |\nabla _{\theta } u|^2  \, h \, d\Theta .
\end{eqnarray*}
By (\ref{intid2}) and (\ref{intid3}) this implies 
\begin{equation}
\label{Jprime}
J'(0)  = 0, 
\end{equation}
and
\begin{equation}
\label{Jprimeprime}
J^{\prime \prime } (0)
 =(k+N+\alpha-1) (k-l-1) \int_{\mathbb{S}_+^{N-1} }  (u+s_1 )^2 \, h \, d\Theta + 
 \int_{\mathbb{S}_+^{N-1} }  |\nabla _{\theta } u|^2 \, h \, d\Theta .
\end{equation}

Now assume that (\ref{isop1}) holds. 
Then we have ${\mathcal R}_{k,l,N,\alpha} (U(t)) \geq {\mathcal R}_{k,l,N,\alpha} (B_{1}^{+} )$ 
for all $t$ with $|t|<t_0 $. In view of (\ref{intid1}) 
this means that $J(t) \geq J(0) $ for $|t|<t_0 $, that is, 
\begin{equation}
\label{Jderiv}
J^{\prime \prime } (0) \geq 0 = J'(0).
\end{equation}
The second condition is (\ref{Jprime}), and the first condition 
implies, in view of (\ref{intid2}) and 
(\ref{Jprimeprime}), that
\begin{eqnarray}
\label{intineq1}
0 & \leq & (k+N+\alpha-1)(k-l-1) \int_{\mathbb{S}_+^{N-1} }   v^2 \, h \,d\Theta + 
 \int_{\mathbb{S} _+^{N-1} } |\nabla _{\theta } v|^2 \, h \, d\Theta 
\\
\nonumber
 & & \forall v\in C^2 (\mathbb{S} _+^{N-1} ) \ \mbox{ with } \ 
 \int_{\mathbb{S}_+^{N-1} } v \, h \, d\Theta =0.
\end{eqnarray}
Applying Proposition 2.1 in \cite{BCM2}, we get
$$
\int_{\mathbb{S} _+^{N-1} } |\nabla _{\theta } v|^2 \, h \, d\Theta \ge (N+\alpha-1)
\int_{\mathbb{S}_+^{N-1} }     v^2  \, h \, d\Theta 
$$
for any $ v\in C^2 (\mathbb{S} _+^{N-1} ) $ with
$ \dint_{\mathbb{S}_+^{N-1} }h v \, d\Theta =0$.
The conclusion follows.
$\hfill \Box $

\section{The case of negative $\alpha$}

\noindent
In this section we firstly show that the relative isoperimetric
problem in $\mathbb{R}_{+}^{2}$ for $\alpha \in \left( -1,0\right) $ and  $k=l=0$ has
no solution. Nevertheless,   in  Theorem \ref{St_Not_Iso}, we prove that,  the second variation
of the perimeter w.r.t. volume-preserving smooth perturbations 
at the half circle is nonnegative for such values 
 of the parameters.

\noindent Throughout this section the points in $\mathbb{R}_{+}^{2}$
will be simply denoted by $(x,y)$.

\bigskip

\begin{theorem}
\label{Not_Ex} 
Let
\begin{equation}
N=2,\text{ }\alpha \in \left( -1,0\right) \,\, \text{and } k=l=0 .
 \label{H_NE}
\end{equation}
Then there is no constant $C\in \left( 0,+\infty \right) $ such that 
\begin{equation*}
\int_{\partial \Omega \backslash \left\{ y=0\right\} }y^{\alpha }dl\geq
C\left( \displaystyle\int\limits_{\Omega }y^{\alpha }dxdy\right) ^{\frac{
\alpha +1}{\alpha +2}},\text{ for any  set }\Omega \subset \mathbb{R}_{+}^{2}.
\end{equation*}
\end{theorem}
\noindent {\sl Proof:} 
\ \ Let $0<a<b$ and 
\begin{equation*}
\Omega _{a,b}:=\left\{ (x,y)\in \mathbb{R}_{+}^{2}:0<x<1,\text{ }a<y<b\right\} .
\end{equation*}
We have 
\begin{equation*}
A_{\alpha }\left( \Omega _{a,b}\right) :=\displaystyle\int\limits_{\Omega
_{a,b}}y^{\alpha }dxdy=\int_{a}^{b}t^{\alpha }dt=\frac{b^{\alpha
+1}-a^{\alpha +1}}{\alpha +1}.
\end{equation*}
while 
\begin{equation*}
P_{\alpha }\left( \Omega _{a,b}\right) :=\int_{\partial \Omega
_{a,b}}y^{\alpha }dl=2\int_{a}^{b}t^{\alpha }dt+a^{\alpha }+b^{\alpha }=
\frac{2}{\alpha +1}\left( b^{\alpha +1}-a^{\alpha +1}\right) +a^{\alpha
}+b^{\alpha }.
\end{equation*}
Setting 
\begin{equation*}
U:=a^{\alpha +1},\text{ \ }V:=b^{\alpha +1}-a^{\alpha +1}\hspace{0.5cm}
(U,V>0)
\end{equation*}
we have 
\begin{equation*}
A_{\alpha }\left( \Omega _{a,b}\right) =\frac{V}{\alpha +1}\ \text{\ and \ }
P_{\alpha }\left( \Omega _{a,b}\right) =\frac{2}{\alpha +1}V+U^{\frac{\alpha 
}{\alpha +1}}+\left( U+V\right) ^{\frac{\alpha }{\alpha +1}}.
\end{equation*}
In order to conclude to proof we claim that $\forall \epsilon >0$ $\exists $ 
$0<a<b$ such that

\begin{equation*}
R_{\alpha }\left( \Omega _{a,b}\right) 
\equiv 
\frac{P_{\alpha }\left(
\Omega _{a,b}\right) }{\left[ A_{\alpha }\left( \Omega _{a,b}\right) \right]
^{\frac{\alpha +1}{\alpha +2}}}<\epsilon .
\end{equation*}
First choose $V$ small enough to have 
\begin{equation*}
2\left( \alpha +1\right) ^{-\frac{1}{\alpha +1}}\text{ }V^{\frac{1}{\alpha +2
}}<\frac{\epsilon }{2}
\end{equation*}
and then $U$ large enough to have 
\begin{equation*}
\frac{U^{\frac{\alpha }{\alpha +1}}+(U+V)^{\frac{\alpha }{\alpha +1}}}{
\left( \frac{1}{\alpha +1}\right) ^{\frac{\alpha +1}{\alpha +2}}V^{\frac{
\alpha +1}{\alpha +2}}}<\frac{\epsilon }{2}.
\end{equation*}
Then 
\begin{equation*}
R_{\alpha }\left( \Omega _{a,b}\right) =2\left( \alpha +1\right) ^{-\frac{1}{
\alpha +1}}\text{ }V^{\frac{1}{\alpha +2}}+\frac{U^{\frac{\alpha }{\alpha +1}
}+(U+V)^{\frac{\alpha }{\alpha +1}}}{\left( \frac{1}{\alpha +1}\right) ^{
\frac{\alpha +1}{\alpha +2}}V^{\frac{\alpha +1}{\alpha +2}}}<\frac{\epsilon 
}{2}+\frac{\epsilon }{2}=\epsilon .
\end{equation*}

\hfill $\Box $

\bigskip 

\noindent Now let $\alpha \in \left( -1,0\right) $ and 
consider the measure $d\nu =\cos ^{\alpha }t \, dt $.
We introduce the weighted
Sobolev space 
$H^{1}\left( \left( -\frac{\pi }{2},\frac{\pi }{2}\right)
;d\nu \right) $ which is made of  functions $\phi :\left( -\frac{\pi }{2}
,\frac{\pi }{2}\right) \rightarrow  \mathbb{R}$ such that
\begin{eqnarray*}
\left\Vert \phi \right\Vert _{H^{1}\left( \left( -\frac{\pi }{2},\frac{\pi }{
2}\right) ; \,d\nu \right) }^{2} 
&=&
\left\Vert \phi \right\Vert_{L^{2}\left(  \left( -\frac{\pi }{2},\frac{\pi }{2} \right) ; \, d\nu \right)
}^{2}+\left\Vert \phi ^{\prime }\right\Vert _{L^{2}\left( \left( -\frac{\pi }{2},
\frac{\pi }{2} \right) ;  \, d\nu \right) }^{2} \\
&=&
\int_{-\frac{\pi }{2}}^{\frac{\pi }{2}}
\phi   (t)^{2} \, d\nu +
\int_{-\frac{\pi }{2}}^{\frac{\pi }{2}}
\phi ^{\prime }(t)^{2} \, d\nu  <\infty .
\end{eqnarray*}
Finally let
\begin{equation*}
V :=\left\{ \phi \in H^{1}\left( \left( -\frac{\pi }{2},\frac{\pi }{2}
\right) ;   \, d\nu    \right) :\int_{-\frac{\pi }{2}}^{\frac{\pi }{2}
}\phi  \, d\nu  =0\right\} .
\end{equation*}
 In the following Lemma we prove that
$ V $
is compactly embedded in 
$L^{2}\left(  \left( -\frac{\pi }{2},\frac{\pi }{2} \right) ;  \, d\nu   \right) $.

\begin{lemma}
\label{embedd}
 If $\left\{ w_{n}\right\} _{n\in
N}\subset V $ is such that
\begin{equation*}
\int_{-\frac{\pi }{2}}^{\frac{\pi }{2}}w_{n}^{\prime }(t)^{2} \, d\nu \leq C\text{ \ }\forall n \in \mathbb{N}
\end{equation*}
then there exists $ w\in V $ such that there holds
\begin{equation*}
\lim_{n\rightarrow \infty }\int_{-\frac{\pi }{2}}^{\frac{\pi }{2}}\left\vert
w_{n}(t)-w(t)\right\vert ^{2}\, d\nu =0\text{.}
\end{equation*}
\end{lemma}

\noindent {\sl Proof:}  \ Note that
\begin{equation*}
\int_{-\frac{\pi }{2}}^{\frac{\pi }{2}}w_{n}^{\prime }(t)^{2}dt\leq \int_{-
\frac{\pi }{2}}^{\frac{\pi }{2}}w_{n}^{\prime }(t)^{2}\cos ^{\alpha }tdt\leq
C\text{ \ }\forall n \in \mathbb{N} .
\end{equation*}
By the definition of $V$ we can infer that for each $n\in \mathbb{N}$, there exists $
t_{n}\in (-\frac{\pi }{2},\frac{\pi }{2})$  such that, up to a subsequence,  $w_{n}(t_{n})=0.$ So
we have
\begin{equation*}
w_{n}(t)=\int_{t_{n}}^{t}w_{n}^{\prime }(\sigma )d\sigma
\end{equation*}
and therefore
\begin{equation*}
\left\vert w_{n}(t)\right\vert ^{2}\leq \left( \int_{-
\frac{\pi }{2}}^{\frac{\pi }{2}}\left\vert w_{n}^{\prime }(\sigma
)\right\vert d\sigma \right) ^{2}\leq \pi \int_{-\frac{\pi }{2}}^{\frac{\pi 
}{2}}\left\vert w_{n}^{\prime }(\sigma )\right\vert ^{2}d\sigma \leq C\text{
\ }\forall n \in \mathbb{N}.
\end{equation*}
So $w_{n}$ is bounded in $H^{1} \left( -\frac{\pi }{2},\frac{\pi }{2}  \right)$
 and, therefore, there exists 
$w\in C^{0}\left( \left[ -\frac{\pi }{2},\frac{\pi }{2
}\right] \right)
\cap H^{1} \left( -\frac{\pi }{2},\frac{\pi }{2}  \right) $ 
such that, up to a subsequence,
\begin{equation*}
w_{n}(t)\rightarrow w(t)\text{ uniformly in }\left[ -\frac{\pi }{2},\frac{
\pi }{2}\right] .
\end{equation*}
The assertion easily  follows, since
\begin{equation*}
\cos ^{\alpha }t\in L^{1}\left( -\frac{\pi }{2},\frac{\pi }{2}\right) \text{
\ }\forall \alpha \in (-1,0) .
\end{equation*}
\hfill $\Box $

\bigskip

\noindent Now define the Rayleigh quotient
\begin{equation*}
Q(v):=\frac{\dint_{-\frac{\pi }{2}}^{\frac{\pi }{2}}v^{\prime }(t)^{2}\cos
^{\alpha }tdt}{\dint_{-\frac{\pi }{2}}^{\frac{\pi }{2}}v(t)^{2}\cos ^{\alpha
}tdt},   \,\,\,\ \text{with} \,\,\,\, v \in V .
\end{equation*}

\begin{lemma} 
\label{W_Wirt}
 There holds  
\begin{equation*}
\mu :=\min_{\phi \in V}Q(v)=1+\alpha .
\end{equation*}
\end{lemma} 

\noindent {\sl Proof:}  \  \ Note that $\sin t\in V $. An integration by parts
gives 
\begin{equation}
Q(\sin t)=\frac{\dint_{-\frac{\pi }{2}}^{\frac{\pi }{2}}\cos ^{\alpha +2}tdt}{
\dint_{-\frac{\pi }{2}}^{\frac{\pi }{2}}\sin ^{2}t\cos ^{\alpha }tdt}
=
\frac{ \left( \alpha +1\right) \dint_{-\frac{\pi }{2}}^{\frac{\pi }{2}}\sin
^{2}t\cos ^{\alpha }tdt}{\dint_{-\frac{\pi }{2}}^{\frac{\pi }{2}}\sin
^{2}t\cos ^{\alpha }tdt}=\alpha +1,  
\label{sint}
\end{equation}
and, therefore 
\begin{equation*}
\mu \leq \alpha +1.
\end{equation*}
Now, by contradiction, assume that 
\begin{equation*}
\mu <1+\alpha .
\end{equation*}
By Lemma \ref{embedd} there exists a function $u\in V$ such that $Q(u)=\mu $
which satisfies  the Euler equation  
\begin{equation}
-\left( u^{\prime }\cos ^{\alpha }(t)\right) ^{\prime }=\mu u\cos ^{\alpha
}(t)\text{ \ on \ }\left( -\frac{\pi }{2},\frac{\pi }{2}\right) .
\label{eig_eq}
\end{equation}
We set 
\begin{equation*}
R(v):=\int_{-\frac{\pi }{2}}^{\frac{\pi }{2}}v^{\prime }(t)^{2}\, d \nu -\mu \int_{-\frac{\pi }{2}}^{\frac{\pi }{2}}v(t)^{2} \, d \nu  ,
\text{ \ }v\in V,
\end{equation*}
and 
\begin{equation*}
u_{1}(t)=\frac{u(t)-u(-t)}{2},\text{ \ }u_{2}(t)=\frac{u(t)+u(-t)}{2}.
\end{equation*}
We have 
\begin{equation*}
R(u)=R(u_{1})+R(u_{2})=0.
\end{equation*}
Hence at least one of the following statements must be true 
\begin{equation}
R(u_{1})\leq 0,  \tag{i}  \label{i}
\end{equation}
or
\begin{equation}  \label{ii}
R(u_{2})\leq 0.  \tag{ii}
\end{equation}

\noindent Our aim is to reach a contradiction by showing that (\ref{i}) 
and (\ref{ii}) are both false.

\vspace{.5cm}

\noindent \textbf{Case (i)}:  Assume $R(u_{1})\leq 0.$

\noindent Since $u_{1}$ is odd we have 
\begin{equation*}
v_{1}:=\frac{u_{1}(t)}{\sin t}\in C^{1}\left( \left[ -\frac{\pi }{2},\frac{
\pi }{2}\right] \right)
\end{equation*}
and 
\begin{equation*}
R(u_{1})=\int_{-\frac{\pi }{2}}^{\frac{\pi }{2}}\left( v_{1}^{\prime }\sin
t+v_{1}\cos t\right) ^{2}\cos ^{\alpha }tdt-\mu \int_{-\frac{\pi }{2}}^{ 
\frac{\pi }{2}}v_{1}^{2}\sin ^{2}t\cos ^{\alpha }tdt=
\end{equation*}
\begin{equation*}
=\int_{-\frac{\pi }{2}}^{\frac{\pi }{2}}2v_{1}^{\prime }v_{1}\sin t\cos
^{\alpha +1}tdt+\int_{-\frac{\pi }{2}}^{\frac{\pi }{2}}\left( v_{1}^{\prime
}\right) ^{2}\sin ^{2}t\cos ^{\alpha }tdt +
\end{equation*}
\begin{equation*}
+\int_{-\frac{\pi }{2}}^{\frac{\pi 
}{2}}v_{1}^{2}\cos ^{\alpha +2}tdt+-\mu \int_{-\frac{\pi }{2}}^{\frac{\pi }{
2 }}v_{1}^{2}\sin ^{2}t\cos ^{\alpha }tdt
\end{equation*}
\begin{eqnarray*}
&=&
(\alpha +1)\int_{-\frac{\pi }{2}}^{ 
\frac{\pi }{2}}v_{1}^{2}\sin ^{2}t\cos ^{\alpha }tdt-\int_{-\frac{\pi }{2}
}^{ \frac{\pi }{2}}v_{1}^{2}\cos ^{\alpha +2}tdt+ \\
&&\int_{-\frac{\pi }{2}}^{\frac{\pi }{2}}\left( v_{1}^{\prime }\right)
^{2}\sin ^{2}t\cos ^{\alpha }tdt+\int_{-\frac{\pi }{2}}^{\frac{\pi }{2}
}v_{1}^{2}\cos ^{\alpha +2}tdt-\mu \int_{-\frac{\pi }{2}}^{\frac{\pi }{2}
}v_{1}^{2}\sin ^{2}t\cos ^{\alpha }tdt
\end{eqnarray*}
Recalling  the assumption  $\alpha +1-\mu >0$, we have 
\begin{eqnarray*}
R(u_{1}) &=&(\alpha +1)\int_{-\frac{\pi }{2}}^{\frac{\pi }{2}}v_{1}^{2}\sin
^{2}t\cos ^{\alpha }tdt+\int_{-\frac{\pi }{2}}^{\frac{\pi }{2}}\left(
v_{1}^{\prime }\right) ^{2}\sin ^{2}t\cos ^{\alpha }tdt-\mu \int_{-\frac{\pi 
}{2}}^{\frac{\pi }{2}}v_{1}^{2}\sin ^{2}t\cos ^{\alpha }tdt \\
&=&(\alpha +1-\mu )\int_{-\frac{\pi }{2}}^{\frac{\pi }{2}}v_{1}^{2}\sin
^{2}t\cos ^{\alpha }tdt+\int_{-\frac{\pi }{2}}^{\frac{\pi }{2}}\left(
v_{1}^{\prime }\right) ^{2}\sin ^{2}t\cos ^{\alpha }tdt\geq 0,
\end{eqnarray*}
where equality holds if and only if $\ \mu =\alpha +1$ and $v_{1}$ is a
constant. This contradicts our assumption.

\vspace{.5cm}

\noindent \textbf{Case (ii)}: Assume $R(u_{2})\leq 0.$

\noindent Since $u_{2}$ is even function belonging to $V$, we have 
\begin{equation*}
0 = \int_{-\frac{\pi }{2}}^{\frac{\pi }{2}}u_{2}\cos ^{\alpha }tdt = 2
\int_{0}^{\frac{\pi }{2}}u_{2}\cos ^{\alpha }tdt.
\end{equation*}
Then there exists $c\in \left( 0,\frac{\pi }{2}\right) $ such that 
\begin{equation*}
u_{2}(c)=u_{2}(-c)=0.
\end{equation*}
From \eqref{eig_eq} we deduce that 
\begin{equation}
\int_{-c}^{c}\left( u_{2}^{\prime }\right) ^{2}\cos ^{\alpha }tdt=
-\int_{-c}^{c}u_{2}\left( u_{2}^{\prime }\cos
^{\alpha }t\right) ^{\prime }dt=\mu \int_{-c}^{c}u_{2}^{2}\cos ^{\alpha }tdt.
\label{-c_+c}
\end{equation}
On the other hand, setting 
\begin{equation*}
v_{2}:=u_{2}\cos ^{\frac{\alpha }{2}}t ,
\end{equation*}
we obtain from (\ref{-c_+c}) 
\begin{eqnarray}
\int_{-c}^{c}\left( u_{2}^{\prime }\right) ^{2}\cos ^{\alpha }tdt
&=&\int_{-c}^{c}\left( v_{2}^{\prime }\cos ^{-\frac{\alpha }{2}}t+\frac{
\alpha }{2}v_{2}\cos ^{-\frac{\alpha }{2}-1}t\sin t\right) ^{2}\cos ^{\alpha
}tdt  \label{v2} \\
&=&\int_{-c}^{c}\left( v_{2}^{\prime }\right) ^{2}dt+\alpha
\int_{-c}^{c}v_{2}v_{2}^{\prime }\tan tdt+\frac{\alpha ^{2}}{4}
\int_{-c}^{c}v_{2}^{2}\tan ^{2}tdt.  \notag
\end{eqnarray}
Since $v_{2}\left( \pm c\right) =0$
and
 $v_{2}\in C^{1}\left[ -c,c\right] $, the
classical one-dimensional Wirtinger inequality implies that 
\begin{equation}
\int_{-c}^{c}\left( v_{2}^{\prime }\right) ^{2}dt\geq \left( \frac{\pi }{2c}
\right) ^{2}\int_{-c}^{c}v_{2}^{2}dt,
 \label{W_1d}
\end{equation}
where equality holds if and only if $v_{2}$ 
is proportional to  $ \sin \left( \dfrac{\pi t}{2c} \right) $

Inequalities (\ref{-c_+c}) and (\ref{W_1d}) ensure 

\begin{eqnarray}
\int_{-c}^{c}\left( u_{2}^{\prime }\right) ^{2}\cos ^{\alpha }tdt &\geq
&\left( \frac{\pi }{2c}\right) ^{2}\int_{-c}^{c}v_{2}^{2}dt \\
&&-\frac{\alpha }{2}\int_{-c}^{c}v_{2}^{2}\left( 1+\tan ^{2}t\right) dt+
\frac{\alpha ^{2}}{4}\int_{-c}^{c}v_{2}^{2}\tan ^{2}tdt  \notag \\
&=&\left( \frac{\pi ^{2}}{4c^{2}}-\frac{\alpha }{2}\right)
\int_{-c}^{c}v_{2}^{2}dt+\left( \frac{\alpha ^{2}}{4}-\frac{\alpha }{2}
\right) \int_{-c}^{c}v_{2}^{2}\tan ^{2}tdt  \notag \\
&>&\left( \frac{\pi ^{2}}{4c^{2}}-\frac{\alpha }{2}\right)
\int_{-c}^{c}v_{2}^{2}dt  \notag \\
&=&\left( \frac{\pi ^{2}}{4c^{2}}-\frac{\alpha }{2}\right)
\int_{-c}^{c}u_{2}^{2}\cos ^{\alpha }tdt.  \notag
\end{eqnarray}
Finally equation \eqref{eig_eq} implies 
\begin{equation*}
1+\alpha >\mu >\frac{\pi ^{2}}{4c^{2}}-\frac{\alpha }{2}\geq 1-\frac{\alpha 
}{2}
\end{equation*}
and therefore \ $\frac{3}{2}\alpha >0 ,$ a contradiction.
\hfill $\Box $

\begin{theorem}
\label{St_Not_Iso}
Let $N=2,\text{ }\alpha \in \left( -1,0\right) $  and  $k=l=0$. Then 
the functional $J$ defined in \eqref{perim},
satisfies $J^{\prime \prime}(0) \geq 0$. 
\end{theorem}

\noindent {\sl Proof:} 
The assertion follows from  Lemma  \ref{W_Wirt} and taking into account of \eqref{Jprimeprime}.
\hfill $\Box $

\vspace{.5cm}

\section{Main results}

This section is devoted to the proof of Theorem \ref{maintheorem}, that is, 
we obtain sufficient conditions on $k,l$ and $N$ such that 
$ C_{k,l,N, \alpha} = C_{k,l,N, , \alpha} ^{rad}$ holds, or equivalently,
\begin{equation}
\label{ineqrad}
{\mathcal R}_{k,l,N, \alpha} (M) \geq C_{k,l,N, \alpha}^{rad}  
\quad \mbox{for all measurable sets $M \subset \mathbb R^{N}_{+}$ with $0< \mu _{l, \alpha} (M) <+\infty $.}
\end{equation}
Proofs of Theorem \ref{ineqrad} are given in various subsections, each of which addresses 
one of the cases ofTheorem \ref{maintheorem}.

First let us recall that the proof of case (i) of Theorem \ref{maintheorem} has been given 
in \cite{ABCMP_atti}.
\begin{remark}
\label{sufficiency}
Condition (\ref{k_l_ineq1}), i.e. $l\frac{N+\alpha-1}{N+\alpha}\le k$ is a 
necessary and sufficient condition for 
$C_{k,l,N, \alpha} >0$.
\end{remark}
{\sl Proof:\/} 
The necessity follows from Lemma \ref{R3}, and the sufficiency in the case 
$l+1\leq k$ follows from case (i) in Theorem \ref{maintheorem}.
Finally, assume that $k< l+1$. Then  (\ref{isopproblem}) is equivalent to 
(\ref{ineqQR}), by Lemma \ref{R2}.
Now the main Theorem of \cite{CKN} tells us that condition 
(\ref{k_l_ineq1}) is also sufficient  
for $C_{k,l,N, \alpha} >0$. 
$\hfill \Box $

\subsection{Proof of Theorem \ref{maintheorem}, case (ii).} 
The case $k \leq 0$ and $\alpha = 0$ has been addressed in  \cite{ChiHo}, Theorem 1.3. 
We significantly extend such a result by considering all nonnegative values of $\alpha$ and treating, 
at least for some values of the parameters,
the equality case in (\ref{isop1}).   
\begin{theorem} 
\label{th1bis} 
Let $k,l$ satisfy 
\begin{equation} \label{lk} 
l \frac{N+\alpha-1}{N+\alpha} \leq k 
\leq \min\{0, l+1\}. 
\end{equation}
Then (\ref{isop1}) holds. 
Moreover if
$l \frac{N+\alpha-1}{N+\alpha} < k$
and 
\begin{equation}
\label{M=BR}
{\mathcal R}_{k,l,N, \alpha} (M) = C_{k,l,N, \alpha} ^{rad} \ 
\mbox{ for some measurable set $M$ with $0<\mu _l (M)< +\infty $},
\end{equation}  
 then $M= B_{R}^{+}$ for some $R>0$. 
\end{theorem}
{\sl Proof  :\/} 
Let $u\in C^{\infty }_0(\mathbb{R}_+ ^N)\setminus \{ 0 \} $.
We set
$$
y:=x|x|^\frac{k}{N+\alpha-1}\, , \quad v(y):=u(x)\, , \quad 
s:=r^\frac{k+N+\alpha-1}{N+\alpha -1}\,
.
$$
Using $N$-dimensional spherical coordinates,  denoting with $\nabla_\theta$ 
the tangential part of the gradient on  
 ${\mathbb S^{N-1}}$, we obtain
\begin{eqnarray} 
\label{cambio1} 
 & &
\int_{\mathbb{R}_+ ^N} x_N^\alpha
|x|^l
|u|^{(l+N+\alpha)/(k+N+\alpha-1)} \, dx
\\
\nonumber
  & = & 
\int_{\mathbb{S}^{N-1}_+} 
\int_0^{\infty} 
r^{l+N+\alpha-1} |u|^{(l+N+\alpha)/(k+N+\alpha-1) }\, 
 h  dr\, d\Theta 
\\
\nonumber
 & = & 
 \frac{N+\alpha-1}{k+N+\alpha-1} 
 \int_{\mathbb{S}^{N-1}_+} 
 \int_0^{\infty} 
 s^{\frac{l+N+\alpha}{k+N+\alpha-1}(N+\alpha-1)-1 }
 |v|^{(l+N+\alpha)/(k+N+\alpha-1)}\, 
   h ds \, d\Theta 
 \\
\nonumber
 & = & 
 \frac{N+\alpha-1}{k+N+\alpha-1} 
 \int_{\mathbb{R} _+^N} y_n^\alpha
 |y|^{\frac{l+N+\alpha}{k+N+\alpha-1}(N+\alpha-1)-N}|v|^{(l+N+\alpha)/(k+N+\alpha-1)}\, dy  
 \\
\nonumber
 & = & 
 \frac{N+\alpha-1}{k+N+\alpha-1} 
 \int_{\mathbb{R} _+^N} 
 |y|^{(l(N+\alpha-1)-k(N+\alpha))/(k+N+\alpha-1)}
 |v|^{(l+N+\alpha)/(k+N+\alpha-1)}\, dy
 \, .
\end{eqnarray}
Further we calculate
\begin{eqnarray}
\label{cambio2}
& &
 \int_{\mathbb{R} _+^N} x_N^\alpha |x|^k |\nabla_x u| \, dx
 \\
  \nonumber
 & = &
 \int_{\mathbb{S}^{N-1} _+} 
 \int_0^{\infty} 
 r^{k+N+\alpha-1}
 \left(  
 u_r ^2 +\frac{|\nabla_\theta  u|^2}{r^2}   
 \right) ^{1/2}h  \, 
 dr \, d\Theta
 \\
 \nonumber
 & = &
 \int_{\mathbb{S}^{N-1} } 
 \int_0^{\infty} 
 s^{N+\alpha-1}
 \left(  
 v_s ^2+\frac{|\nabla_\theta v|^2}{s^2} 
 \left(
 \frac{N+\alpha-1}{k+N+\alpha-1} \right) ^2  \right) ^{1/2} \,  h \, ds \, d\Theta
\nonumber 
\\
\nonumber
 & \geq &
 \int_{\mathbb{S}^{N-1} } 
 \int_0^{\infty} 
 s^{N+\alpha-1}
 \left(  
 v_s ^2 +\frac{|\nabla_\theta v|^2}{s^2}  \right) ^{1/2} \, h \, ds \, d\Theta
\\
\nonumber 
 & = & 
 \int_{\mathbb{R}_+^N}y_N^\alpha |\nabla_y v| \, dy \, ,
\end{eqnarray}
where we have used (\ref{lk}).
By \eqref{cambio1} and \eqref{cambio2} we deduce,
 \begin{eqnarray}
 \label{Q2}
  & & \hspace {1cm}
 {\mathcal Q}_{k,l,N, \alpha}(u)
 \\
 \nonumber
  & \geq & 
 \frac{\displaystyle 
 \int_{\R ^N_+} y_N^\alpha |\nabla_y v| \, dy}{\displaystyle 
 \left( 
 \int_{\R ^N_+} y_N^\alpha    |y|^{l' }|v|^{(l+N+\alpha)/(k+N+\alpha-1)}\, dy  \right) 
 ^{(k+N+\alpha-1)/(l+N+\alpha)} 
 }
 \left(
\frac{k+N+\alpha-1}{N+\alpha-1}
\right) ^{(k+N+\alpha-1)/(l+N+\alpha)}
\\
\nonumber
 & = &
\left( \frac{k+N+\alpha-1}{N+\alpha-1} \right)
^{(k+N+\alpha-1)/(l+N+\alpha)} 
{\mathcal Q}_{0,l' ,N, \alpha }(v)\, ,
\end{eqnarray}
where we have set $l' :=\frac{l(N+\alpha-1)-k(N+\alpha)}{k+N+\alpha-1}$. 
Note that we have $-1 \leq l' \leq 0$ by the assumptions (\ref{lk}).
\\
Hence we may apply Lemma \ref{R2} to both sides of (\ref{Q2}). 
This yields
\begin{equation}
\label{relationCC}
C_{k,l,N, \alpha} \geq \left( \frac{k+N+\alpha-1}{N+\alpha-1} \right) ^{(k+N+\alpha-1)/(l+N+\alpha)} C_{0,l', N, \alpha} .
\end{equation}
Furthermore, Lemma \ref{rangekl1} tells us that
\begin{equation}
\label{CCrad}
C_{0,l',N, \alpha} = C_{0,l', N, \alpha} ^{rad} .
\end{equation}
Since also 
$$
\left( \frac{k+N+\alpha-1}{N+\alpha-1} \right)
^{(k+N+\alpha-1)/(l+N+\alpha)} C_{0,l',N, \alpha} ^{rad} =C_{k,l,N, \alpha}^{rad} \, .
$$
From this, (\ref{relationCC}) and (\ref{CCrad}), 
we deduce that  $C_{k,l,N, \alpha}\ge C_{k,l,N, \alpha}^{rad}$. 
Since $C_{k,l,N, \alpha}\le C_{k,l,N, \alpha}^{rad}$ by definition, (\ref{isop1}) follows.
\\
Next assume that
${\mathcal{R}}_{k,l,N, \alpha} (M) = C_{k,l,N, \alpha} ^{rad}$ for
some measurable set $M \subset \mathbb R^{N}_{+}$ with $0<\mu _l (M)<
+\infty $.  
If $l(N+\alpha-1)/(N+\alpha) <k$, then Lemma \ref{rangekl1} tells us that we must have 
$M=B_{R}^{+} $  for some $R>0$.
$\hfill \Box $

\begin{remark} 
\rm
$\text{}$
\\
 {\bf (a)} A well-known special case of Theorem \ref{th1bis} is $k=0 = l $, 
see   \cite{MadernaSalsa},  \cite{BCM} and \cite{XR}.
\\ 
{\bf (b)}
The idea to use spherical coordinates, 
and in particular the inequality (\ref{cambio2}) in our last proof, 
appeared already in some work of T. Horiuchi, 
see \cite{H} and \cite{HK}. 
\end{remark}

\subsection{Proof of Theorem \ref{maintheorem}, case (iii).} 
Now we treat the case when $k$ assumes non-negative values. 
Throughout this subsection we assume  $k\leq l+1$. 
The main result is Theorem \ref{th1ter}. 
Its proof is long and  requires some auxiliary results. 
But the crucial idea is an interpolation argument that occurs 
in the proof of the following Lemma \ref{4.3}, formula (\ref{ineq2}). 
\begin{lemma}  
\label{4.3}
Assume $l(N+\alpha-1)/(N+\alpha)\leq k$ and $k\geq 0$.   
Let $u\in C_0 ^1 (\mathbb{R}^N_+)\setminus \{ 0 \} $, $u\geq 0$, 
and define $y,z$ and $v$ by   
\begin{equation}
\label{transf1}
 y:= x|x| ^{\frac{k}{N+\alpha-1}} , \ z:= |y| \ \mbox{ and }\  v(y) := 
  u(x),  \qquad   x\in \mathbb{R}_+ ^N .
\end{equation}
Then for every
$A\in \left[ 
0, \frac{(N+\alpha-1) ^2}{ (k+N+\alpha-1 )^2 } 
\right] 
$,
\begin{equation}
\label{ineq1}
 {\mathcal Q}_{k,l,N, \alpha} (u)  \geq \left( 
\frac{k+N+\alpha-1}{N+\alpha-1} 
\right) 
^{ \frac{k+N+\alpha-1}{l+N+\alpha} } 
 \cdot
\frac{ 
\left( 
\dint_{\mathbb{R}_+ ^N } 
y_N^\alpha|\nabla _y v| \, 
\, dy 
\right) 
^A
\cdot 
\left( 
\dint_{\mathbb{R}_+ ^N } 
y_N^\alpha | v_z|  
\, dy  
\right) 
^{1-A}
 }{ 
\left( 
\dint_{\mathbb{R}_+ ^N } 
y_N^\alpha |y|
^{ \frac{l(N+\alpha-1)-k(N+\alpha)}{k+N+\alpha-1} } 
v 
^{ \frac{l+N+\alpha}{k+N+\alpha-1} } 
\, dy 
\right) 
^{ \frac{k+N+\alpha-1}{l+N+\alpha} }
} 
.
\end{equation}
\end{lemma}    

{\sl Proof:} 
We calculate as in the proof of Theorem \ref{th1bis} ,
$$
\int_{\mathbb{R}_+^N }x_N^\alpha |x| ^k |\nabla _x u | \, dx 
 =    
 \int_{\mathbb{S}_+^{N-1} } 
 \int_0^{\infty} 
 s^{N+\alpha-1}
 \left(  
 v_s ^2+\frac{|\nabla_\theta v|^2}{s^2} 
 \left(
 \frac{N+\alpha-1}{k+N+\alpha-1} \right) ^2  \right) ^{1/2} \, h \, ds \, d\Theta
  $$
 Since the mapping
 $$
 t\longmapsto \log 
 \left(
 \int _{{\mathbb S}_+^{N-1} } 
 \int_0^{+\infty } z^{N+\alpha-1}  
 \sqrt{ v_z ^2 + t\frac{|\nabla _{\theta } v|^2 }{z^2 }  } \, h \, dz\, d\Theta \right)
 $$
 is concave, we deduce that for every 
 $A\in \left[ 0, \frac{(N+\alpha-1) ^2}{ (k+N+\alpha-1 )^2 } \right] $,
\begin{eqnarray}
\label{ineq2}
 & & \int_{\mathbb{R}_+^N }x_N^\alpha |x| ^k |\nabla _x u | \, dx 
\\ 
\nonumber
 & \geq & 
 \left( 
 \int_{\mathbb{S}_+^{N-1} } \int_0 ^{+\infty }  z^{N+\alpha-1}  
 \sqrt{ 
 v_z ^2 + \frac{|\nabla _{\theta } v|^2 }{z^2 } 
  } 
\, h  \, dz \, d\Theta \right)
  ^A 
  \cdot 
\left(  
\int_{\mathbb{S}_+^{N-1}}
\int_0^{+\infty } z^{N+\alpha-1}  |v_z | \, h  \, dz\, d\Theta 
\right) 
^{1-A}
\\
\nonumber
 & = & 
 \left( 
 \int_{\mathbb{R}^N_+ } y_N^\alpha|\nabla _y v| \, dy 
 \right) 
 ^A
  \cdot 
 \left( 
 \int_{\mathbb{R}^N_+ } y_N^\alpha|v_z | \, dy 
 \right) 
 ^{1-A} .
 \end{eqnarray}
 Finally, we have 
 \begin{equation}
 \label{equaldenom}
 \int_{\mathbb{R}^N_+ }x_N^\alpha |x| ^l u
 ^{ \frac{l+N+\alpha}{k+N+\alpha-1} } \, dx   = 
 \frac{N+\alpha -1}{k+N+\alpha -1} \int_{\mathbb{R}_+^N} y_N^\alpha
 |y| 
 ^{ \frac{l(N+\alpha -1)-k(N+\alpha )}{k+N+\alpha -1} } 
 v
 ^{ \frac{l+N+\alpha }{k+N+\alpha -1} } \, dy .
 \end{equation}
Now (\ref{ineq1}) follows from 
(\ref{ineq2}) and (\ref{equaldenom}).
$\hfill \Box $
\\[0.1cm]   
Next we want to estimate the right-hand-side of (\ref{ineq1}) from below. 
We will need a few more properties of the starshaped rearrangement.
\begin{lemma}  
\label{4.2}  
Assume $l(N+\alpha-1)/(N+\alpha)\leq k$. 
Then we have for any function 
$v\in C_0 ^1 (\mathbb{R}^N_+ )\setminus \{ 0 \}$ with $v\geq 0$, 
\begin{eqnarray}
\label{starsh5}
 & & 
 \int_{\mathbb{R}_+^N } y_N^\alpha |y| 
 ^{ \frac{l(N+\alpha-1)-k(N+\alpha)}{k+N+\alpha-1} } 
 v
 ^{ \frac{l+N+\alpha}{k+N+\alpha-1} } 
 \, dy 
\leq 
\int_{\mathbb{R}_+^N } y_N^\alpha |y| 
^{ \frac{l(N+\alpha-1)-k(N+\alpha)}{k+N+\alpha-1} } 
\widetilde{v}
^{ \frac{l+N+\alpha}{k+N+\alpha-1} } 
\, dy,
\\
\label{vL1}
 & & 
\frac{ 
y\cdot \nabla \widetilde{v} }{|y|} 
\equiv 
\frac{
\partial \widetilde{v} }{
\partial z } 
\in L^1 (\mathbb{R}_+ ^N )
\quad \mbox{and }
\\
\label{starsh6}
 & &
\int_{\mathbb{R}^N_+ } 
y_N^\alpha\left| \frac{ \partial v}{\partial z}  \right| \, dy 
\geq  
\int_{\mathbb{R}^N_+ } 
 y_N^\alpha\left|\frac{ \partial \widetilde{v} }{\partial z} \right|
\, dy.
\end{eqnarray}
\end{lemma}

{\sl Proof:} Let us prove (\ref{starsh5}). Set 
$$
w(y):= |y| 
^{ \frac{l(N+\alpha-1)-k(N+\alpha)}{l+N+\alpha} } . 
$$
Since $l(N+\alpha-1)-k(N+\alpha)\leq 0$, we have $w= \widetilde{w} $. 
Hence (\ref{starsh5}) follows from (\ref{harlit}) and (\ref{monrearr}). 
\\
Next let $\zeta := z^N  $ and define $V$ and $\hat{V} $ by 
$V(\zeta ,\theta) := v(z\theta )$, and 
$\widehat{V} (\zeta ,\theta ) := \widetilde{v} (z\theta )$. 
Observe that for each $\theta \in \mathbb{S}_+^{N-1} $,  
$\widehat{V} (\cdot , \theta ) $ is the equimeasurable 
non-increasing rearrangement of 
$V (\cdot ,\theta )$. Further we have 
$$
\frac{ \partial v }{ \partial z }
  =  
 N\zeta 
 ^{ \frac{N-1}{N}  } 
 \frac{ \partial V}{\partial \zeta } \ \mbox{ and } \ 
\frac{ \partial \widetilde{v} }{ \partial z }
  = 
 N\zeta 
 ^{ \frac{N-1}{N}  } 
 \frac{ \partial \widehat{V} }{ \partial \zeta  } 
 .
$$
Since $\frac{\partial v}{\partial z} \in L^{\infty } (\mathbb{R}^N )$, 
Lemma \ref{Landes} tells us 
that for every $\theta \in {\mathbb S} ^{N-1} $,
\begin{eqnarray*} 
\int_0^{+\infty } z^{N+\alpha-1} \left| \frac{\partial v}{\partial z} 
(z\theta ) \right| \, dz 
 & = & 
 \int_0 ^{+\infty } \zeta 
 ^{ \frac{N+\alpha-1}{N} } 
 \left| \frac{\partial V}{\partial \zeta } (\zeta ,\theta )
 \right| \, d \zeta  
 \\
 & \geq & 
 \int_0 ^{+\infty } \zeta 
 ^{ \frac{N+\alpha-1}{N} } 
 \left| \frac{\partial \widehat{V} }{\partial \zeta }  
 (\zeta ,\theta )\right| \, d\zeta 
 \\
  & = & 
\int_0^{+\infty } z^{N+\alpha-1} \left|
\frac{ \partial \widetilde{v} }{\partial z} (z\theta )\right| \, dz .
\end{eqnarray*}
Integrating this over ${\mathbb S}_+ ^{N-1}$, 
we obtain (\ref{starsh6}).
$\hfill \Box$

A final ingredient is
\begin{lemma}  
\label{4.1}
Assume that $l(N+\alpha-1)/(N+\alpha)\leq k$, 
and let $M \subset \mathbb R^{N}_{+}$ be a bounded starshaped set. Then
\begin{eqnarray}
\label{holder1}
 & & 
 \left( 
 \int_M y_N^\alpha |y| 
 ^{ \frac{l(N+\alpha-1) -k(N+\alpha)}{k+N+\alpha-1} } 
 \, dy 
 \right) 
 ^{ \frac{k+N+\alpha-1}{l+N+\alpha} } 
 \\
 \nonumber
  & \leq & 
  d_1 
  \left( 
  \int_M y_N^\alpha\, dy 
  \right) 
  ^{  \frac{(N+\alpha-1)(l-k+1) }{ l+N+\alpha}   } 
  \cdot 
  \left( 
  \int_M y_N^\alpha|y|^{-1} \, dy 
  \right) 
  ^{  \frac{k(N+\alpha)-l(N+\alpha-1) }{ l+N+\alpha}  } ,
\quad \mbox{ where}
\\
 & & 
\label{d1}
d_1 = 
\left( 
\frac{k+N+\alpha-1}{l+N+\alpha} 
\right) 
^{ \frac{k+N+\alpha-1}{l+N+\alpha} } 
\cdot 
\left( 
\frac{N+\alpha}{N+\alpha-1} 
\right) 
^{ \frac{(N+\alpha-1)(l-k+1)}{l+N+\alpha} } .
\end{eqnarray}
Moreover,
if $k<l+1$ and $l(N+\alpha-1)/(N+\alpha) <k$, then 
equality in (\ref{holder1}) holds only
if $M=B_{R}^{+} $ for some $R>0$.
\end{lemma}    

{\sl Proof:} Since $M$ is starshaped, 
there is a bounded measurable function 
$m : \mathbb{S} ^{N-1}_{+} \to [0, +\infty )$, such that
\begin{equation} 
\label{Mrepresent}
M= \{ z\theta :\, 0\leq z < m (\theta ), \ 
\theta \in \mathbb{S} ^{N-1}_{+} \} .
\end{equation} 
Using H\"older's inequality we obtain
\begin{eqnarray}
\label{chain}
 & & 
 \hspace{1cm}  \int_M y_N^\alpha|y| 
 ^{ \frac{l(N+\alpha-1) -k(N+\alpha)}{k+N+\alpha-1} } 
 \, dy 
 \\
 \nonumber
 & = & 
 \frac{k+N+\alpha-1}{(l+N+\alpha)(N+\alpha-1)} 
\int_{{\mathbb S}_+ ^{N-1} } m (\theta ) 
^{ \frac{(l+N+\alpha)(N+\alpha-1)}{k+N+\alpha-1} } 
\, h \, d\Theta 
\\
 \nonumber
 & = & 
 \frac{k+N+\alpha-1}{(l+N+\alpha)(N+\alpha-1)} 
\int_{{\mathbb S}_+ ^{N-1} } m (\theta ) 
^{ \frac{k(N+\alpha)-l(N+\alpha-1)}{k+N+\alpha-1}(N+\alpha-1) } m (\theta ) 
^{ \frac{(N+\alpha-1)(l-k+1)}{k+N+\alpha-1} (N+\alpha)} 
\, h \, d\Theta 
\\
\nonumber
& \leq &    
 \frac{k+N+\alpha-1}{(l+N+\alpha)(N+\alpha-1)} 
 \left( 
\int_{{\mathbb S}_+ ^{N-1} } 
m (\theta ) ^{N+\alpha} 
\, h \, d\Theta
\right) 
^{ \frac{(N+\alpha-1)(l-k+1)}{k+N+\alpha-1} } 
\\
\nonumber
 \qquad & \times
 &  
\left( 
\int_{{\mathbb S}_+ ^{N-1} } 
m (\theta ) ^{N+\alpha-1} 
\, h \, d\Theta
\right) 
^{ \frac{k(N+\alpha)- l(N+\alpha-1)}{k+N+\alpha-1} }
\\
\nonumber
 & = &  
 \frac{k+N+\alpha-1}{(l+N+\alpha)(N+\alpha-1)} 
 \left( 
 (N+\alpha) \int_M y_N^\alpha dy 
 \right) 
 ^{ \frac{(N+\alpha-1)(l-k+1)}{k+N+\alpha-1} } 
  \times
 \\
\nonumber
 \qquad & \times
 &  
 \left( 
 (N+\alpha-1) \int_M |y| ^{-1} y_N^\alpha\, dy 
 \right) 
 ^{ \frac{k(N+\alpha)- l(N+\alpha-1)}{k+N+\alpha-1} } ,
\end{eqnarray}
and (\ref{holder1}) follows.
If $k<l+1$ and $l(N+\alpha -1)/(N+\alpha ) < k$, then (\ref{chain}) holds   
with equality only if $m (\theta )=\mbox{const }$.   
$\hfill \Box $
\\[0.1cm] 

Now we are ready to prove our main result. 
\begin{theorem}   
\label{th1ter}
Assume $0\le k\leq l+1$ and
\begin{equation}
\label{crucial}
l\leq 
\frac{(k+N+\alpha-1)^3 }{(k+N+\alpha-1)^2 - \frac{(N+\alpha-1)^2 }{ N+\alpha}} -N -\alpha.
\end{equation}
Then (\ref{isop1}) holds.
Furthermore, if 
 inequality (\ref{crucial}) is strict,
then (\ref{M=BR}) holds only if $M=B_{R}^{+}$ for some $R>0$.  
\end{theorem}    

{\sl Proof: } First observe that the conditions 
$k\geq 0$ and (\ref{crucial}) also imply
$l(N+\alpha-1)/(N+\alpha) \leq k$. Let $u \in C_0 ^{\infty } 
(\mathbb{R}_+^N)\setminus \{ 0\} $, $u\geq 0$,  
and let $v$ be given by 
(\ref{transf1}).
In view of (\ref{crucial}), we may choose
$$
A=\frac{(N+\alpha)(l-k+1)}{l+N+\alpha} 
$$
to obtain 

\begin{eqnarray}
\label{ineq5bis}
 {\mathcal Q}_{k,l,N, \alpha} (u) 
& \geq &
\left(  \frac{k+N+\alpha-1}{N+\alpha-1}   \right)^{ \frac{k+N+\alpha-1 }{ l+N+\alpha} } \times
 \\
 \nonumber
 &  \times &
\frac{ 
\left( 
\dint_{\mathbb{R} ^N } 
y_N^\alpha |\nabla _y v|  
\, dy 
\right) 
^{ \frac{(N+\alpha)(l-k+1) }{ l+N+\alpha} } 
\cdot 
\left( 
\dint_{\mathbb{R}_+ ^N } 
y_N^\alpha |v_z |
\, dy  
\right) 
^{ \frac{k(N+\alpha)-l(N+\alpha-1) }{ l+N+\alpha} } 
 }{
\left( 
\dint_{\mathbb{R}_+^N } y_N^\alpha |y|
^{ \frac{l(N+\alpha-1)-k(N+\alpha ) }{ k+N+\alpha-1} }
v 
^{ \frac{l+N+\alpha}{k+N+\alpha-1} } 
\, dy 
\right) 
^{ \frac{k+N+\alpha-1}{l+N+\alpha} } 
 }
\end{eqnarray} .

Further, (\ref{starsh6}) and 
Hardy's inequality yield
\begin{equation}
\label{ineq4}
 \int_{\mathbb{R}^N } y_N^\alpha |v_z| \, dy
  \geq   
  \int_{\mathbb{R}^{N}_{+} } 
  y_N^\alpha |\widetilde{v}_z| \, dy
 \geq  (N+\alpha -1) 
  \int_{\mathbb{R}^{N}_{+} } y_N^\alpha\frac{\widetilde{v}}{|y|}  
  \, dy\,,
 \end{equation}
 where $\widetilde{v}$ denotes the starshaped rearrangement of $v$.
Together with (\ref{ineq5bis}) and (\ref{starsh5}) this leads to
\begin{eqnarray}
\label{ineq5final}
 {\mathcal Q}_{k,l,N, \alpha} (u) 
 & \geq & 
(N+\alpha-1)
^{ \frac{k(N+\alpha)-l(N+\alpha-1)}{l+N+\alpha} } 
 \left( 
 \frac{k+N+\alpha-1}{N+\alpha-1} 
 \right) 
^{ \frac{k+N+\alpha-1 }{ l+N+\alpha} } 
\cdot 
\\
\nonumber
 & & \cdot
\frac{ 
\left( 
\dint_{\mathbb{R}_+ ^N } 
 y_N^\alpha|\nabla _y v|  
\, dy 
\right) 
^{ \frac{(N+\alpha)(l-k+1) }{ l+N+\alpha} } 
\cdot 
\left( 
\dint_{\mathbb{R}_+^N } 
 y_N^\alpha\frac{\widetilde{v} }{|y|}  
\, dy  
\right) 
^{ \frac{k(N+\alpha)-l(N+\alpha-1) }{ l+N+\alpha} } 
 }{
\left( 
\dint_{\mathbb{R}_+^N } y_N^\alpha |y|
^{ \frac{l(N+\alpha-1)-k(N+\alpha) }{ k+N+\alpha-1} }
\widetilde{v} 
^{ \frac{l+N+\alpha}{k+N+\alpha-1} } 
\, dy 
\right) 
^{ \frac{k+N+\alpha-1}{l+N+\alpha} } 
 }  .
\end{eqnarray} 
Now let $M $ be a bounded measurable subset of $\mathbb R^{N}_{+}$.
Then combining (\ref{limperim}), (\ref{limmeas}) 
and the argument leading to (\ref{CklNsmooth}) we deduce that  
there exists a sequence of non-negative functions 
$\{ u_n \} \subset C_0 ^1 (\mathbb{R}^N_+  )$ such that
\begin{equation}
\label{lim1}
\lim_{n\to \infty } 
\int_{\mathbb{R}_+^N } x_N^\alpha|x| ^k |\nabla u_n | \, dx = 
P_{\mu _k, \alpha } (M) 
\end{equation}
and
\begin{equation}
\label{lim2}
 u_n \longrightarrow \chi_{M } \quad 
\mbox{ in $L^p (\mathbb{R}^{N}_{+} ) $ for every $p\geq 1 $.}
\end{equation}
We define  
$M ':= \{ y= x|x|^{\frac{k}{N+\alpha-1}} : \, x \in M \} $ 
and $v_n (y) := u_n (x) $.

Let $\widetilde{v_n } $ and $\widetilde{M '} $ 
be the starshaped rearrangements of $v_n $ and $M ' $ respectively. 
Then (\ref{lim1}) and (\ref{lim2}) also imply 
\begin{eqnarray}
\label{lim3}
 & & \lim_{n\to \infty } \int_{\mathbb{R}^N_+ } y_N^\alpha |\nabla _y v_n | \, dy  = 
P_{\mu _0, \alpha } (M'  ),
\quad \mbox{and}
\\ 
\label{lim4}
 & & \widetilde{v_n }\longrightarrow 
\chi _{\widetilde{M ' } } \
\mbox{ in $L^p (\mathbb{R}^{N}_{+} ) $ for every $p\geq 1 $.}
\end{eqnarray}
Choosing $u=u_n $ in (\ref{ineq5final}) and passing to the limit 
$n\to \infty $, 
we obtain, using (\ref{lim1}), (\ref{lim2}), (\ref{lim3}), 
(\ref{lim4}) and Proposition  \ref{BCM2} 
\begin{eqnarray}
\label{ineq6}
 & & {\mathcal R}_{k,l,N, \alpha} (M )
 \\
 \nonumber
 & \geq &  
(N+\alpha-1)
^{ \frac{k(N+\alpha)-l(N+\alpha-1)}{l+N+\alpha} } 
\left( 
 \frac{k+N+\alpha-1}{N+\alpha-1} 
\right) 
^{ \frac{k+N+\alpha-1}{l+N+\alpha} } 
\cdot 
\\
\nonumber
 & & \cdot
\frac{ 
\left( P_{\mu _0, \alpha } (\widetilde{M'}) \right)
^{ \frac{(N+\alpha)(l-k+1)}{l+N+\alpha} } 
\cdot 
\left( 
\dint_{\widetilde{M' } }  
\frac{y_N^\alpha dy}{|y| } 
\right) 
^{ \frac{k(N+\alpha)-l(N+\alpha-1)}{l+N+\alpha} } 
}{
\left( 
\dint_{\widetilde{M '}  } y_N^\alpha |y|
^{ \frac{l(N+\alpha-1)-k(N+\alpha)}{k+N+\alpha-1} } 
\, dy 
\right) 
^{ \frac{k+N+\alpha-1}{l+N+\alpha} } 
} 
\\
 \nonumber
 & \geq &  
(N+\alpha-1)
^{ \frac{k(N+\alpha)-l(N+\alpha-1)}{l+N+\alpha} } 
\left( 
 \frac{k+N+\alpha-1}{N+\alpha-1} 
\right) 
^{ \frac{k+N+\alpha-1}{l+N+\alpha} }
\left( 
C_{0,0,N, \alpha} ^{rad}
\right) 
^{ \frac{(N+\alpha)(l-k+1)}{l+N+\alpha} }  
\times 
\\
\nonumber
 & & \times 
\frac{ 
\left( \mu_{0,\alpha} (\widetilde{M'}) \right)
^{ \frac{(N+\alpha-1)(l-k+1)}{l+N+\alpha} } 
\cdot 
\left( 
\dint_{\widetilde{M' } }  
\frac{ y_N^\alpha dy}{|y| } 
\right) 
^{ \frac{k(N+\alpha)-l(N+\alpha-1)}{l+N+\alpha} } 
}{
\left( 
\dint_{\widetilde{M '}  } y_N^\alpha |y|
^{ \frac{l(N+\alpha-1)-k(N+\alpha)}{k+N+\alpha-1} } 
\, dy 
\right) 
^{ \frac{k+N+\alpha-1}{l+N+\alpha} } 
} .
\end{eqnarray}
In view of (\ref{holder1}) and since 
$\mu _0 (M') = \mu_0 (\widetilde{M'})$
we finally get from this
\begin{eqnarray}
\label{ineq7}
& &
{\mathcal R}_{k,l,N, \alpha} (M) 
 \\
 &\geq&
  \nonumber
(N+\alpha-1)
^{ \frac{k(N+\alpha)-l(N+\alpha-1)}{l+N+\alpha} } 
\left( \frac{k+N+\alpha-1}{N+\alpha-1} \right) 
^{ \frac{k+N+\alpha-1}{l+N+\alpha} } 
\left( 
C_{0,0,N, \alpha} ^{rad}
\right) 
^{ \frac{(N+\alpha)(l-k+1)}{l+N+\alpha} }  \frac{1}{d_1} 
\\
\nonumber 
 & = &  \left( \dint_{\mathbb S^{N-1}_{+}}  \,h \,d \Theta \right) ^{\frac{l-k+1}{l+N+\alpha} } 
\cdot (l+N+\alpha)^{\frac{k+N+\alpha-1}{l+N+\alpha} } = C_{k,l,N, \alpha} ^{rad}  ,
\end{eqnarray}
and (\ref{isop1}) follows by (\ref{CklNsmooth}).
\\
Now assume that (\ref{M=BR}) holds. 
If inequality (\ref{crucial}) is strict, 
then Lemma \ref{rangekl1} tells us that we must have 
$M= B_{R}^{+} $ for some $R>0$.
$\hfill \Box $
\begin{remark}  
\rm 
Note that if $N+\alpha \geq 3$
, then (\ref{crucial}) covers the important range 
$$
l=0\leq k\leq 1. 
$$
However, we emphasize that this is not true when $2\leq N+\alpha <3$.  
\end{remark}    
\medskip

\noindent

\section{Applications}
In this section we provide some applications of our results.
\subsection{P\'{o}lya-Szeg\"o principle}
First we obtain a P\'{o}lya-Szeg\"o principle related to our 
isoperimetric inequality (\ref{isop1}) (cf. \cite{Talenti2})
Assume that the numbers $k, l$ and $\alpha$ satisfy (\ref{ass1}) and one of the conditions {\bf (i)}-{\bf (iii)} of 
Theorem \ref{maintheorem}. Then (\ref{mainineq}) implies
\begin{equation} 
\label{Isop_klalpha}
\int_{\partial \Omega }
|x|^k x_N ^{\alpha } {\mathcal H}_{N-1}(dx)
\geq 
\int_{\partial \Omega ^{ \star }}
|x|^k x_N ^{\alpha } 
{\mathcal H}
_{N-1}(dx)
\end{equation}
for every smooth set $\Omega \subset \mathbb{R} ^N_+ $, where $\Omega ^{\star}$ is the $\mu_{l,\alpha } $-symmetrization of $\Omega $.
We will use (\ref{Isop_klalpha}) to prove the following
\begin{theorem} 
\label{ps}
(P\'{o}lya-Szeg\"o principle) 
Let the numbers $k,l$ and $\alpha $ satisfy one of the conditions {\bf (i)}-{\bf (iii)} 
of Theorem \ref{maintheorem}.
Further, let $p\in [1, +\infty)$ and $m:= pk+(1-p)  l $. Then there holds 
\begin{equation}
\int_{\mathbb{R}^N _+ }
\left\vert  \nabla u\right\vert ^p 
d \mu _{m,\alpha } (x)
\geq 
\int_{\mathbb{R}^N _+ }\left\vert \nabla
u^{ \star }\right\vert ^{p}
d\mu_{m,\alpha } (x)
\quad 
\forall u\in {\mathcal D} ^{1,p} (\mathbb{R} ^N _+ , d\mu _{m, \alpha } ),
\label{PS_k_l}
\end{equation}
where $u^{\star } $ denotes the $\mu _{l,\alpha } $-symmetrization of $u$.
\end{theorem}
{\sl Proof:} 
It is sufficient to consider the case that $u$ is non-negative. Further,
by an approximation argument we may assume that
$u \in C^{\infty}_{0}(\mathbb{R}^{N} ) $.
Let 
\begin{eqnarray*}
I & := &  
\int_{\mathbb{R}^{N}_+ } | \nabla u| ^{p}
|x| ^{pk+(1-p)l} x_N ^{\alpha }\, dx \quad \mbox{and}\\
I ^{\star }  & := &  
\int_{\mathbb{R}^{N}_+ } | \nabla u^{\star} | ^{p}
|x| ^{pk+(1-p)l} x_N ^{\alpha }\, dx .
\end{eqnarray*}
The coarea formula yields 
\begin{eqnarray}
\label{1coarea}
I 
 & = & 
 \int_{0}^{\infty }\int_{u=t} |\nabla
u| ^{p-1} |x| ^{pk+(1-p)l} x_N ^{\alpha }\, {\mathcal H}_{N-1}(dx)\, dt \quad \mbox{and}
\\
\label{coarea2}
I^{\star}
 & = & 
 \int_{0}^{\infty }\int_{u^{\star } =t} |\nabla
u^{\star} | ^{p-1} |x| ^{pk+(1-p)l} x_N ^{\alpha }\, {\mathcal H}_{N-1}(dx)\, dt 
.
\end{eqnarray}
Further, H\"older's inequality gives
\begin{equation}
\label{1holder}
\int_{ u=t} |x|^k x_N ^{\alpha } \, {\mathcal H} _{N-1} (dx)
\leq 
\left( \int_{ u=t} |x|^{kp +l(1-p)} |\nabla u| ^{p-1} x_N ^{\alpha } \, {\mathcal H} _{N-1} (dx) \right) ^{\frac{1}{p} } 
\cdot
\left( \int_{ u=t} \frac{|x|^l x_N ^{\alpha }}{|\nabla u|} \, {\mathcal H}_{N-1} (dx) 
\right)
^{\frac{p-1}{p} } ,
\end{equation} 
for a.e. $t\in [0, +\infty )$. 
Hence (\ref{1coarea}) together with (\ref{1holder}) tells us that
\begin{equation}
\label{coarea3}
I 
\geq 
\int_{0}^{\infty }
\left( \int_{u=t} |x| ^{k} x_N ^{\alpha }\, {\mathcal H}
_{N-1}(dx)
\right) ^{p} \cdot 
\left( 
\int_{u=t}\frac{ |x| ^{l}x_N ^{\alpha }}{
| \nabla u| } x_N ^{\alpha } \, {\mathcal H}_{N-1}(dx)
\right) ^{1-p} \, dt.
\end{equation}
Since $u^{\star} $ is a radial function, we obtain in an analogous manner,
\begin{equation}
\label{coarea4}
I^{\star} 
=
\int_{0}^{\infty }
\left( \int_{u^{\star} =t} |x| ^{k} x_N ^{\alpha } \, {\mathcal H}
_{N-1}(dx)
\right) ^{p} \cdot 
\left( 
\int_{u^{\star} =t}\frac{ |x| ^{l}x_N ^{\alpha } }{
| \nabla u^{\star} | } \, {\mathcal H}_{N-1}(dx)
\right) ^{1-p} \, dt.
\end{equation}
Observing that
\begin{equation}
\label{meas_u>t}
\int_{u>t} |x|^{l} x_N ^{\alpha } \, dx
=
\int_{u^{\star }>t}
|x|^{l} x_N ^{\alpha } \, dx \quad  \forall t\in [0, +\infty ), 
\end{equation}
Fleming-Rishel's formula yields
\begin{equation}
\label{flemingrishel} 
\int_{u=t } \frac{|x|^l x_N ^{\alpha }}{|\nabla u|} \, {\mathcal H}_{N-1} (dx) 
=
\int_{u^{\star} =t } \frac{|x|^l x_N ^{\alpha }}{|\nabla u^{\star} |} \, {\mathcal H}_{N-1} (dx)
\end{equation}
for a.e. $t\in [0, +\infty )$. 
Hence
(\ref{flemingrishel}) and (\ref{Isop_klalpha}) give
\begin{eqnarray*}
 & &
\int_{0}^{\infty }
\left( \int_{u=t} |x|^k x_N ^{\alpha } \, {\mathcal H}
_{N-1}(dx) \right) ^{p}
\cdot 
\left( \int_{u=t}\frac{| x| ^{l} x_N ^{\alpha } }{
| \nabla u| } \, {\mathcal H}_{N-1}(dx) \right) ^{1-p}
\, dt 
\\
 & \geq &
\int_{0}^{\infty }\left( \int_{u^{\star} =t} |x| ^{k} x_N ^{\alpha }
\, {\mathcal H}_{N-1}(dx) \right) ^{p} \cdot \left( \int_{u^{\star}=t}
\frac{|x|^{l} x_N ^{\alpha } }{|\nabla u^{\star} | } \, {\mathcal H}
_{N-1}(dx)\right) ^{1-p} \, dt.
\end{eqnarray*}
Now (\ref{PS_k_l})  follows from this, (\ref{coarea3}) and (\ref{coarea4}).
$\hfill \Box$
\\[0.1cm]
An important particular case of Theorem \ref{ps} is 
\begin{corollary}
\label{specialcasePS}
Let $p\in [1, +\infty )$, $N+\alpha \geq 3 $, $a\geq 0 $, $u\in {\mathcal D} ^{1,p} 
(\mathbb{R}^N _+ , d\mu _{ap ,\alpha }) $, and let $u^{\star } $ be the $\mu_{0,\alpha } $-symmetrization of $u$.
Then
\begin{equation}
\label{PSspecial}
\int_{\mathbb{R}^N _+ } \left| \nabla u\right|^p \, d\mu _{ap, \alpha } (x) 
\geq 
\int_{\mathbb{R}^N _+ } \left| \nabla u^{\star} \right|^p \, d\mu _{ap, \alpha } (x) .
\end{equation}
\end{corollary}
{\sl Proof: } We choose $k:= a $ and $l:= 0$. If $a\in [0,1]$ then $k,l$ 
satisfy either one of the conditions {\bf (ii)} or {\bf (iii)}, see also Remark 5.2. If $a\geq 1 $, then $k,l$ satisfy condition {\bf (i)} of Theorem \ref{maintheorem}. Hence (\ref{PSspecial}) follows from Theorem \ref{ps}.
$\hfill \Box $
\subsection{Caffarelli-Kohn-Nirenberg-type inequalities}
Next we will use Theorem \ref{ps} to obtain best constants in some 
inequalities of Caffarelli-Kohn-Nirenberg-type. 

Let $p,q, a, b$ be real numbers  such that
\begin{eqnarray}
\label{CKNassump1}
 & & 1\leq p \leq q \left\{ 
 \begin{array}{ll} 
 \leq \frac{(N+\alpha )p}{N+\alpha -p} &  \mbox{ if } \ p< N+\alpha 
\\
 < +\infty & \mbox{ if } \ p\geq N + \alpha 
 \end{array}
 \right.
 , 
 \\
\nonumber
 & & a> 1-\frac{N+\alpha }{p}, \quad \mbox{and }
 \\
\nonumber
 & & b= b(a,p,q,N, \alpha ) = (N+\alpha ) \left( \frac{1}{p} -\frac{1}{q} \right) + a-1 .
\end{eqnarray}
We define
\begin{eqnarray}
\label{p*}
p^* & := & \left\{ 
\begin{array}{ll} 
\frac{(N+\alpha )p}{N+\alpha -p} & \mbox{ if } p<N+\alpha 
\\
+\infty & \mbox{ if } p\geq N+\alpha 
\end{array}
\right.
,
\\
& &\nonumber
 \\
\label{fctalE}
E_{a,p,q,N, \alpha } (v)
 & := & 
\frac{\dint_{\mathbb{R} ^N _+ } |x|^{ap}  |\nabla v|^p x_N ^{\alpha }\, dx
 }{
\left( \dint_{\mathbb{R}^N _+ } |x|^{bq} |v|^q x_N ^{\alpha } \, dx \right) ^{p/q}  }, 
\quad v\in C_0 ^{\infty } (\mathbb{R}^N  )\setminus \{ 0\}  ,
\\
\label{SapqN}
S_{a,p,q,N, \alpha } & := & \inf \{ E_{a,p,q,N,\alpha } (v): \, v\in C_0 ^{\infty } 
(\mathbb{R}^N  ) \setminus \{ 0\} \}, \quad \mbox{and}
\\
\label{SapqNrad}
S_{a,p,q,N,\alpha } ^{rad} & := & \inf \{ E_{a,p,q,N,\alpha } (v): \, v\in C_0 ^{\infty } 
(\mathbb{R}^N  )\setminus \{ 0\} ,  \ v \mbox{ radial }\}.
\end{eqnarray}
Note that with this new notation we have
\begin{eqnarray}
\label{E=Q}
E_{k,1,\frac{l+N+\alpha }{k+N+\alpha -1} ,N, \alpha } (v) & = & {\mathcal Q}_{k,l,N,\alpha } (v) \quad \forall v\in C_0 ^{\infty } (\mathbb{R}^N  )\setminus \{ 0\} ,
\\
\label{S=C}
S_{k,1,\frac{l+N+\alpha }{k+N+\alpha -1} ,N, \alpha } (v) & = & C_{k,l,N,\alpha } \quad \mbox{and}
\\
\label{Srad=Crad}
S_{k,1,\frac{l+N+\alpha }{k+N+\alpha -1} ,N,\alpha } ^{rad} & = & C_{k,l,N,\alpha } ^{rad} .
\end{eqnarray}
\\
We are interested 
in the range of values $a$ (depending on $p,q,N$ and $\alpha $) for which 
\begin{equation}
\label{S=S_rad} 
S_{a,p,q,N,\alpha } = S_{a,p,q,N,\alpha } ^{rad}
\end{equation} 
holds. 
\\
First observe that the case $1<p=q$ (which is equivalent to $a-b=1$) 
corresponds to a weighted Hardy-Sobolev-type inequality. Note that inequality \eqref{eq:theorem:Hardy with weight} below was already known when $\alpha=0$ (see, for example \cite{HK} and references therein). We have:
\begin{theorem}
\label{hardysobolev}
\label{theorem:Hardy with weight}
Let $p\geq 1$, $\alpha\geq 0$ and $k\in\mathbb{R}$ be such that  $N-p+\alpha +k>0$. 
Then we have 
\begin{equation}
 \label{eq:theorem:Hardy with weight}
  \int_{\mathbb{R}^N_+} |\nabla u(x)|^p \, d\mu _{k,\alpha } (x) 
\geq
\left(\frac{N-p+k+\alpha  }{p}\right)^p
\int_{\mathbb{R}^N_+ } \frac{| u(x)|^p }{|x|^p } \, d\mu_{k,\alpha } (x)
\end{equation}
for all $u\in {\mathcal D} ^{1,p} (\mathbb{R}^N_+ , d\mu_{k,\alpha }) $
and
\begin{equation}
\label{constant}
S_{a,p,p,N,\alpha } ^{rad} = S_{a,p,p,N,\alpha } 
  =\left(\frac{N-p+k+\alpha  }{p}\right)^p .
\end{equation}
Moreover there is no function $u\in {\mathcal D} ^{1,p}(\mathbb{R}^N_+,d\mu_{k,\alpha } )$ satisfying equality in \eqref{eq:theorem:Hardy with weight} and such that\\
$\int _{\mathbb{R}^N_+ } |\nabla u|^p d\mu_{k,\alpha } \neq 0.$
\end{theorem}

{\sl Proof:} The first two steps follow the line of proof of \cite{GaPe}, Lemma 2.1.
\\
\textit{Step 1.} Assume first that $u\in C_0^{\infty}(\mathbb{R}^N)$. Then
we have for every $x\in \mathbb{R}^N _+ $,
$$
  |u(x)|^p=
- \int_1^{\infty}\frac{d}{dt}|u(tx)|^p\, dt=
  - \int_1^{\infty} p|u(tx)|^{p-2}u(tx)\langle x,\nabla u(tx)\rangle \, dt .
$$
Multiplying this with $x_N ^{\alpha } |x|^{k-p} $ and integrating over $\mathbb{R}^N _+$ we find
\begin{eqnarray}
\nonumber
 \int_{\mathbb{R}^N_+ }|u(x)|^p  x_N ^{\alpha } |x|^{k-p} \, dx & = & - p\int_{1}^{\infty}\left[
 \int_{\mathbb{R}^N_+ } |u(tx)|^{p-2} u(tx) \langle x, \nabla u(tx)\rangle x_N ^{\alpha } |x|^k \, dx
 \right] \, dt
\\
 \nonumber
& = &
 - p\int_{1}^{\infty}\frac{1}{t^{N-p+\alpha +k }}\left[
 \int_{\mathbb{R}^N_+ } \frac{|u(y)|^{p-2} u(y) }{|y|^{p}}\langle y, \nabla u(y)\rangle y_N ^{\alpha } |y|^k \, dy
 \right] \, dt 
 \\
\label{identityhardy}
& =& 
- \frac{p}{N-p+\alpha +k }
 \int_{\mathbb{R}^N_+} \frac{|u(x)|^{p-2} u(x) }{|x|^{p}}\langle x, \nabla u(x)\rangle x_N ^{\alpha } |x|^k \, dx .
\end{eqnarray}
Note that by a density argument (\ref{identityhardy}) still holds for functions $u\in {\mathcal D}^{1,p} (\mathbb{R}^N_+ , d\mu _{k,\alpha } )$. 
In view of the inequality
\begin{equation}
 \label{eq:estimate nabla u by Cauch-Sch}
- u(x) \langle x,\nabla u(x)\rangle \leq |u(x)||x| |\nabla u(x)|
\end{equation}
this leads to  
\begin{equation}
\label{ineq1hardy}
 \int_{\mathbb{R}^N_+ }|u(x)|^p x_N ^{\alpha }|x|^{k-p} \, dx \leq
 \frac{p}{N-p+k+\alpha }
 \int_{\mathbb{R}^N_+ } \frac{|u(x)|^{p-1} }{|x|^{p-1}}|\nabla u(x)| x_N ^{\alpha } |x|^k \, dx .
\end{equation}
Using H\"older's inequality, with $p'$  being the conjugate exponent of $p$, we obtain that (this step is not necessary if $p=1$)
\begin{eqnarray}
\nonumber
 & & \int_{\mathbb{R}^N_+ } \frac{|u(x)|^{p-1}}{|x|^{p-1}}|\nabla u(x)|x_N ^{\alpha } |x|^k \, dx 
\\
\nonumber
 &
= &
 \int_{\mathbb{R}^N_+}\left\{ 
 \frac{|u(x)|^{p-1}}{|x|^{p-1}}\left[ x_N ^{\alpha } |x|^k \right] ^{1/p'} \right\} \left\{ |\nabla u(x)|\left[ x_N ^{\alpha } |x|^k \right] ^{1/p} \right\} \, dx
\\
\label{ineq2hardy}
 & \leq &
 \left(
\int_{\mathbb{R}^N_+ }|u(x)|^p x_N ^{\alpha } |x|^{k-p}\, dx
\right)^{1/p'}
\cdot \left(
\int_{\mathbb{R}^N_+ }|\nabla u(x)|^p x_N ^{\alpha } |x|^k \, dx
\right)^{1 /p} .
\end{eqnarray}
Plugging this estimate into (\ref{ineq1hardy}) concludes the first statement of the theorem.
\smallskip

\textit{Step 2.} Next we show (\ref{constant}).
 Let $\varepsilon >0$ and define 
$$
  M_{\epsilon}=\frac{N-p+k+\alpha +\epsilon}{p},\qquad
  u_{\epsilon}(x)=\left\{\begin{array}{rl}
                          1&\text{ if }|x|\leq 1
                          \smallskip \\
                          |x|^{-M_{\epsilon}}&\text{ if }|x|>1.
                         \end{array}\right.
$$
Note that
$$
  \int_{\mathbb{R}^N_+}|\nabla u_{\epsilon}|^p x_N ^{\alpha } |x|^k \, dx ={M_{\epsilon}}^p\int_{\mathbb{R}^N_+ \backslash B_1}x_N ^{\alpha }|x|^{k-(M_{\epsilon}+1) p}\, dx.
$$
Hence, by Lemma \ref{lemma:integrability w times power} (ii) below we obtain for any $\epsilon >0$ that $u_{\epsilon}\in {\mathcal D}^{1,p}(\mathbb{R}^N_+, d\mu_{k,\alpha } ).$ 
On the other hand, we have that
$$
  \int_{\mathbb{R}^N_+}|u_{\epsilon}(x)|^p x_N ^{\alpha }|x|^{k-p}\, dx=
  \int_{\mathbb{R}^N_+ \backslash B_1} x_N ^{\alpha } |x|^{k-(M_{\epsilon}+1)p}\, dx +\beta,
$$
where, by Lemma \ref{lemma:integrability w times power} (i),
$$
  \beta=\int_{B_1^+ } x_N ^{\alpha }|x|^{k-p}<\infty.
$$
Now set
\begin{displaymath}
  \displaystyle{
  Q_{\epsilon}
=
\frac{
\int_{\mathbb{R}^N_+ } |\nabla u_{\epsilon}|^p x_N ^{\alpha } |x|^{k} \, dx
}{
 \int_{\mathbb{R}^N_+} |u_{\epsilon}|^p x_N ^{\alpha } |x|^{k-p } \, dx }=
  \frac{
\int_{\mathbb{R}^N_+ \backslash B_1} x_N ^{\alpha } |x|^{k- (M_{\epsilon}+1)p}\, dx
}{
\beta +\int_{\mathbb{R}^N_+\backslash B_1}|x|^{k-(M_{\epsilon}+1)p}}
\, dx  .}
\end{displaymath}
Note also that $(M_{\epsilon}+1)p=N+k+\alpha +\epsilon$.
 Therefore we obtain from Lemma \ref{lemma:integrability w times power} (iii) that
$$
  \lim_{\epsilon\to 0}Q_{\epsilon}=(M_0)^p=\left(\frac{N-p+k+\alpha }{p}\right)^p.
$$
This proves the second equality in (\ref{constant}). The first equality in (\ref{constant}) follows from the fact that the approximating functions $u_{\varepsilon}$ are radial.  

\textit{Step 3.} Let us now show that there is no nontrivial function satisfying equality in \eqref{eq:theorem:Hardy with weight}.
\\
Assume that equality holds in (\ref{eq:theorem:Hardy with weight}). Then there holds equality in (\ref{ineq1hardy}) and (\ref{ineq2hardy}). Hence we must have
\begin{eqnarray}
\label{identity3hardy}
 & & -u(x) \langle x, u(x)\rangle =|u(x)||x|\, |\nabla u(x)| \quad \mbox{and} 
\\
\label{identity4hardy}
 & & \frac{|u(x)|}{|x|} = \frac{p}{N-p+k+\alpha } \, |\nabla u(x)| \quad \mbox{for a.e. $x\in \mathbb{R}^N _+ .$}
\end{eqnarray}
An integration of this leads to
\begin{equation}
\label{u=}
u(x) = |x|^{-(N-p+k+\alpha )/p} h\left( x|x|^{-1} \right) ,
\end{equation}
with a measurable function $h: \mathbb{S} ^{N-1} _+ \to \mathbb{R}$.
Since $|x|^{-1} u\in L^p (\mathbb{R}^N _+, d\mu_{k,\alpha }) $, this implies that $h=0$ a.e. on $\mathbb{S}^{N-1} _+ $. The claim is proved.
$\hfill \Box$
\noindent
\begin{lemma}  
\label{lemma:integrability w times power}
Let $\delta >0$. Then
\begin{eqnarray*}
\mbox{(i)} & & 
  \int_{B_1 ^+ }x_N ^{\alpha } |x|^{-N -\alpha +\delta }\, dx <\infty, \quad \mbox{ and }
\\
\mbox{(ii)} & & 
  \int_{\mathbb{R}^N_+ \backslash B_1} x_N ^{\alpha } |x|^{-N -\alpha -\delta }\, dx <\infty.
\end{eqnarray*}
Further, there holds
$$
  \lim_{\delta \to 0+0 }\int_{\mathbb{R}^N_+ \backslash B_1} x_N ^{\alpha }|x|^{-N -\alpha -\delta }\, dx =\infty.
$$
\end{lemma}
{\sl Proof: }
We use $N$-dimensional spherical coordinates to show that
\begin{align*}
  \int_{B_1 ^+ } x_N ^{\alpha } |x|^{-N -\alpha +\delta }=&\int_{\mathbb{S}^{N-1}_+ }\left(
\int_0^1 \left( \frac{x}{|x|}\right) ^{\alpha } r^{-1+\delta } dr \right)d\mathcal{H}^{N-1}(x)
  \smallskip \\
  =&\int_{\mathbb{S}^{N-1}_+ } \left(\frac{x}{|x|}\right) ^{\alpha }d\mathcal{H}^{N-1}(x)\left(\int_{0}^1 r^{-1+\delta }dr\right).
\end{align*}
From this (i) follows. (ii) and (iii) follow similarly.
$\hfill \Box$ 
\\[0.1cm]
\hspace*{1cm} From now on let us assume that
\begin{equation}
\label{maincase}
1<p<q \left\{ 
\begin{array}{ll}
 \leq p^* & \mbox{ if }\ p<N+\alpha 
\\
 <+\infty & \mbox{ if } \ p\geq N+\alpha 
\end{array}
\right.
.
\end{equation} 
We begin with the following
\begin{lemma}
\label{CKN}
Assume that $a, b, p,q,N$ and $ \alpha $ satisfy the conditions (\ref{CKNassump1}) and (\ref{maincase}).  
Further, assume that there exist real numbers $k$ and $l$ which satisfy $l+N+\alpha >0$ and one of the conditions 
{\bf (i)}-{\bf (iii)} of Theorem \ref{maintheorem}, and such that
\begin{eqnarray}
\label{akl}
 & & ap = kp + l(1-p) \ \mbox{ and }
 \\
\label{bq<l}
 & & bq \leq l.
\end{eqnarray} 
Then (\ref{S=S_rad}) holds.
\end{lemma}
{\sl Proof:} Let $u\in {\mathcal D} ^{1,p} (\mathbb{R} ^N _+, d\mu_{ap, \alpha } )\setminus \{ 0\} $,  
and let $u^{\star} $ be the $\mu_{l,\alpha} $-symmetrization of $u$. 
Then we have by Theorem \ref{ps} and (\ref{akl}), 
\begin{equation}
\label{ps1}
\int_{\mathbb{R} ^N _+ } |x|^{ap} |\nabla u| ^p x_N ^{\alpha } \, dx \geq  
\int_{\mathbb{R} ^N _+ } |x|^{ap} |\nabla u^{\star}| ^p x_N ^{\alpha } \, dx.
\end{equation}
Further, it follows from (\ref{hardylitt1}) and (\ref{bq<l})
that 
\begin{equation}
\label{bqint}
\int_{\mathbb{R} ^N _+ } |x|^{bq} | u| ^q x_N ^{\alpha } \, dx \leq  
\int_{\mathbb{R} ^N } |x|^{bq} | u^{\star}| ^q x_N ^{\alpha } \, dx.
\end{equation}
Finally, (\ref{ps1}) together with (\ref{bqint}) yield
\begin{equation}
\label{E>E*}
E_{a,p,q,N,\alpha } (u) \geq E_{a,p,q,N,\alpha } (u^{\star} ),
\end{equation}
and the assertion follows.
$\hfill \Box $
\\[0.1cm] 
Now we define
\begin{eqnarray}
\label{def_a1}
a_1 & := & \frac{N+\alpha -1}{q-\frac{q}{p} +1} +1 -\frac{N+\alpha }{p}, \ \ \mbox{ and }
\\
\label{def_a2}
a_2 & := &  
\frac{N+\alpha -1}{(q- \frac{q}{p} +1 )\sqrt{ (N+\alpha )( \frac{1}{p} -\frac{1}{q})}} +1 -\frac{N+\alpha }{p} .
\end{eqnarray}
Observe that the conditions (\ref{maincase}) imply that 
\begin{equation}
\label{a2>a1>0}
a_2\geq a_1 
\geq 0,
\end{equation}
and equality in the two inequalities holds iff $p<N+\alpha $ and $q=p^* $. 
\\
Moreover, an elementary calculation shows that
\begin{eqnarray}
\label{a1cond}
a_1 & = &
\max \Big\{ a: \, a= k + l\left( \frac{1}{p} -1\right) , \  bq\leq l  ,
\\
\nonumber
 & & \qquad \qquad
-N-\alpha < l \leq k \frac{N+\alpha }{N+\alpha -1 }  \leq 0 \Big\} 
 \quad \mbox{ and }
\\
\label{a2cond}
a_2 & = &
\max \Big\{ a: \, a= k + l\left( \frac{1}{p} -1\right) , \ bq\leq l ,
\ k\geq 0, 
\\
\nonumber
 & & \qquad \qquad 0< l+ N+\alpha \leq  \frac{  (k+N+\alpha -1)^3}{(k+N+\alpha -1)^2   - \frac{(N+\alpha -1)^2 }{N+\alpha } }   \Big\} .
\end{eqnarray} 
The main result of this section is the following
\begin{theorem}
\label{best_a}
Assume that (\ref{maincase}) holds.
Then we have
\begin{equation}
\label{s=s*}
S_{a,p,q,N,\alpha } = S_{a,p,q,N,\alpha } ^{rad} \qquad \forall a\in \Big(
1-\frac{N+\alpha }{p} ,a_2 \Big].
\end{equation}
\end{theorem}
{\sl Proof: }  Let $a \in \Big( 
1-\frac{N+\alpha }{p} ,a_2 \Big]$. We define 
\begin{eqnarray}
\label{l}
l & := & q \left( a+ \frac{N+\alpha }{p} -1 \right) -N -\alpha , \quad 
\mbox{and }
\\
\label{k} 
k & := & \left( 1+ q-\frac{q}{p} \right)  \left( a+ \frac{N+\alpha }{p} -1 \right) -N-\alpha +1 .
\end{eqnarray}
This implies 
\begin{eqnarray*}
a & = & k+l \left(\frac{1}{p} -1 \right) ,  
\\
bq & = & l \quad \mbox{and} 
\\
l+ N+\alpha  & = & \frac{k+ N+\alpha -1 }{ \frac{1}{q} -\frac{1}{p} +1} >0 . 
\end{eqnarray*}
Now we split into two cases:
\\[0.1cm]
{\bf 1.} Let $a\leq a_1 $. 
\\
Then  
$$
k\leq 0,
$$
and since $q\leq p^* $ if $p< N+\alpha $ and $q<+\infty $ otherwise, we have 
\begin{eqnarray*}
l\frac{N+\alpha -1}{N+\alpha } -k 
& = & (k+N+\alpha -1) 
\frac{
-\frac{1}{N+\alpha } -\frac{1}{q} +\frac  {1}{p} }{ \frac{1}{q}-\frac{1}{p} +1} 
\\
 & \leq & 0.
\end{eqnarray*}
Hence we are in case {\bf (ii)} of Theorem \ref{maintheorem}, so that the assertion follows by Lemma \ref{CKN}, for $a \leq a_1 $.
\\[0.1cm]
{\bf 2.} Next let $a_1 \leq a\leq a_2 $.
\\
This implies
\begin{eqnarray}
\nonumber
k & \geq & 0 \quad \mbox{and }
\\
\label{estk}
k+ N+\alpha -1 & \leq & 
\frac{N+\alpha -1}{\sqrt{ (N+\alpha ) \left( \frac{1}{p}-\frac{1}{q} \right) } } .
\end{eqnarray}
Now, from (\ref{estk}) we deduce
\begin{eqnarray*}
 & & l+N+\alpha - \frac{ 
(k+N+\alpha -1 ) ^3 
}{
(k+N+\alpha -1 )^2 - \frac{(N+\alpha -1)^2 }{N+\alpha } 
}
\\
 & = & \frac{ 
(k+N+\alpha -1) 
\left( 
(k+N+\alpha -1)^2 
\left( \frac{1}{p}-\frac{1}{q} 
\right) 
- \frac{(N+\alpha -1)^2 }{N+\alpha } \right) 
}{
\left( 
\frac{1}{q} -\frac{1}{p}+1
\right) 
\left(  
(k+N+\alpha -1)^2- \frac{(N+\alpha -1)^2 }{N+\alpha } 
\right) 
} 
\\
 & \leq & 0.
\end{eqnarray*}
Hence we are in case {\bf (iii)} of Theorem \ref{maintheorem}, so that the assertion follows again by Lemma \ref{CKN} .
$
\hfill \Box$
\\[0.1cm] 
{\bf Remark 6.1:} The characterizations (\ref{a1cond}) and (\ref{a2cond}) and the inequalities (\ref{a2>a1>0}) show that 
 the bound $a_2 $
cannot be improved using our method.
\\[0.1cm]
Finally we evaluate  the constants $S_{a,p,q,N,\alpha } ^{rad} $ 
and  the corresponding radial minimizers.
\\
For any radial function $v\in C_0 ^{\infty } 
(\mathbb{R}^N  ) \setminus \{ 0\} $, it is easy to check the following equality
$$
E_{a,p,q,N, \alpha } (v)
=
\left[B\left( \frac{
N-1}{2},\frac{\alpha +1}{2}\right) 
\right]^{1-\frac pq}
\frac{
\pi ^{\frac{N-1}{2}\frac{q-p}{q}}}{
\left(\Gamma \left[ 
\frac{N-1}{2}\right)\right]^\frac{q-p}{q} }
\frac{\dint_{\mathbb{R} ^N _+ } |x|^{ap+\alpha}  |\nabla v|^p \, dx
 }{
\left( \dint_{\mathbb{R}^N _+ } |x|^{bq+\alpha} |v|^q  \, dx \right) ^{p/q}  }, 
$$ 
Therefore by Theorem 1.4 in \cite{Musina}, we deduce that the function
$$
U(x)=\left(1+|x|^\frac{(N-p+ap+\alpha)(q-p)}{p(p-1)}\right)^\frac{p}{p-q}\,.
$$
achieves the infimum of $E_{a,p,q,N, \alpha }$, that is   $S_{a,p,q,N,\alpha } ^{rad}=E_{a,p,q,N, \alpha } (U)$.

\subsection{Problems in an orthant} Among the possible extensions of our isoperimetric results we would like to address  a problem in an orthant with monomial weights. 
Let $O_+ $ denote the orthant 
$$
O_+ := \{ x\in \mathbb{R} ^N :\, x_i >0 , \, i=1,\ldots , N \} ,
$$
and let $a _1 , \ldots , a _N $ be positive numbers.
Using multi-index notation
we have
\begin{eqnarray*}
{\bf a } & := & (a _1 , \ldots , a _N ),
\\
| {\bf a } | & := & a _1 + \ldots  + a _N ,
\\
x^{{\bf a}}  & := & x_1 ^{a_1 } \cdots x_N ^{a_N } , \quad (x\in \mathbb{R}^N ).
\end{eqnarray*}
Following the lines of proof of Theorem 1.1 we obtain the following isoperimetric result. We leave the details to the reader.
\begin{theorem}
\label{secondmaintheorem}
Let $N\in \mathbb{N} $, $N\geq 2$, 
$k,l \in \mathbb{R} $, ${\bf a} = (a _1 , \ldots , a _N ) $ where $a _i >0 $, ($i=1, \ldots ,N$), and $l+N+|{\bf a }| >0$. 
Further, assume that one of the following conditions holds:
\\
{\bf (i)} $l+1\leq k $;
\\
{\bf (ii)} $k\leq l+1$ and $ l\frac{N+|{\bf a}| -1}{N+|{\bf a}| } \leq k\leq 0$; 
\\  
{\bf (iii)} $N\geq 2$, $ 0\leq k\leq l+1$ and 
\begin{equation}\label{l_1N3new}
l\leq   
\frac{(k+N+|{\bf a}| -1)^3 }{(k+N+|{\bf a}|-1)^2 - \frac{(N+|{\bf a}|-1)^2 }{N+|{\bf a}| } } -N -|{\bf a}| \,.
\end{equation}
\\
Then  
\begin{equation}
\label{mainineqnew}
\dint_{\partial \Omega } |x|^k x^{{\bf a}}\, {\mathcal H}_{N-1} (dx)
\geq 
D   
\left( 
\dint_{\Omega } |x|^l x^{{\bf a}}\, dx 
\right) 
^{(k+N+|{\bf a}|-1)/(l+N+|{\bf a}|) } , 
\end{equation}
for all smooth sets $\Omega $ in 
$O_+   $,
where 
\begin{eqnarray}
\label{defCklnew}
D= D(k,l,N, {\bf a} ) & := &
\frac{\dint_{\partial B_1 } |x|^k x^{{\bf a}}\, {\mathcal H}_{N-1} (dx)}
{\left( \dint_{B_1 \cap  O_+ }
 |x|^l x^{{\bf a}} \, dx \right ) ^{(k+N+|{\bf a}|-1)/(l+N+|{\bf a}|) } } .
\end{eqnarray}
Equality in (\ref{mainineqnew}) holds  if $\Omega =B_R\cap O_+ $.
\end{theorem}

\section*{Acknowledgements} 
The authors are grateful to Gyula Csat\' o who kindly comunicated to us a proof of a general Hardy type inequality a particular case of which is Theorem \ref{hardysobolev}. The authors would to thanks University of Naples Federico II and South China University of Technology of Guangzhou for supporting some visiting appointment and their colleugues for their kind hospitality.


\bigskip



\begin{thebibliography}{99} 
\renewcommand{\baselinestretch}{0.9}
\small 



\bibitem{ABCMP}
{\sc  A. Alvino, F. Brock, F. Chiacchio, A. Mercaldo, M.R. Posteraro},
Some isoperimetric inequalities on $\mathbb R^N$ with respect to weights
$|x|^{\alpha}$, {\sl J. Math. Anal. Appl.} {\bf 451}, no. 1,  (2017), 280--318.

\bibitem{ABCMP_atti}
{\sc  A. Alvino, F. Brock, F. Chiacchio, A. Mercaldo, M.R. Posteraro}, On  weighted isoperimetric inequalities with non-radial densities (2018), ArXiv:1804.02282v1.


\bibitem{AH} 
{\sc  H. Ando, T. Horiuchi},
 On the weighted rearrangement of functions and degenerate nonlinear elliptic equations,
  {\sl Math. J. Ibaraki Univ.} {\bf 44} (2012), 17--31.


\bibitem{BCMR}
{\sc V. Bayle, A. Ca\~{n}ete, F. Morgan, C. Rosales}, On the isoperimetric
problem in Euclidean space with density. {\sl Calc. Var. PDE} {\bf 31} (2008), 27--46.

\bibitem{BBMP} 
{\sc M.F. Betta, F. Brock, A. Mercaldo, M.R. Posteraro}, 
A weighted isoperimetric inequality and applications to symmetrization. 
{\sl J. Inequal. Appl.} {\bf 4} (1999), no. 3, 215--240. 

\bibitem{BBMP2} 
{\sc M.F. Betta, F. Brock, A. Mercaldo, M.R. Posteraro}, 
Weighted isoperimetric inequalities on $\mathbb{R}^N $ and applications to rearrangements.
{\sl Math. Nachr.} {\bf 281} (2008), no. 4, 466--498.

\bibitem{BBCLT}   
{\sc W. Boyer, B. Brown, G. Chambers, A. Loving, S. Tammen}, 
Isoperimetric regions in $\mathbb{R}^n $ with density $r^p $, {\sl Anal. Geom. Metr. Spaces} 
{\bf 4} (2016), 236--265.  

\bibitem{BrasPhil} 
{\sc L. Brasco, G. De Philippis, B. Ruffini},  Spectral optimization for the Stekloff-Laplacian: 
the stability issue. {\sl J. Funct. Anal.} {\bf 262} (2012), no. 11, 4675--4710.

\bibitem{BCM} 
{\sc F. Brock, F. Chiacchio, A. Mercaldo}, A class of
degenerate elliptic equations and a Dido's problem with respect to a
measure. {\sl J. Math. Anal. Appl.} {\bf 348} (2008), no. 1, 356--365.

\bibitem{BCM2} 
{\sc F. Brock, F. Chiacchio, A. Mercaldo}, Weighted
isoperimetric inequalities in cones and applications. {\sl Nonlinear
Analysis T.M.A.} {\bf 75} (2012), no. 15, 5737--5755.

\bibitem{BCM3} 
{\sc F. Brock, F. Chiacchio, A. Mercaldo}, 
A weighted isoperimetric inequality in an orthant. 
{\sl Potential Anal.} {\bf 41} (2012),   171--186.


\bibitem{BMP} 
{\sc F. Brock, A. Mercaldo, M.R. Posteraro}, 
On isoperimetric inequalities with respect to infinite measures.
{\sl Revista Matem\'{a}tica Iberoamericana} {\bf 29} (2013), 665--690.

\bibitem{XR} 
{\sc X. Cabr\' e, X. Ros-Oton}, Sobolev and isoperimetric
inequalities with monomial weights. 
{\sl J. Differential Equations } {\bf 255}  (2013),  4312--4336. 


\bibitem{CKN} 
 {\sc L.  Caffarelli, R.  Kohn, L. Nirenberg},
 First order interpolation inequalities with weights. {\sl Compositio Math.} 
 {\bf 53} (1984), no. 3, 259--275.

\bibitem{CaldMus}
{\sc P. Caldiroli, R. Musina}, Symmetry Breaking of Extremals for the
Caffarelli-Kohn-Nirenberg Inequalities
in a Non-Hilbertian Setting, {\sl Milan J. Math.} {\bf 81} (2013), 421--430.

 
 \bibitem{CMV} 
{\sc A. Ca\~{n}ete, M. Miranda Jr., D. Vittone}, Some
isoperimetric problems in planes with density. {\sl J. Geom. Anal.} 
{\bf 20} (2010), no.2, 243--290.

 
 \bibitem{CJQW} 
{\sc T. Carroll, A. Jacob,  C. Quinn, R. Walters}, The isoperimetric problem on planes with density.
 {\sl Bull. Aust. Math. Soc. } {\bf 78}  (2008), no.2, 177--197.

 
 \bibitem{CRT} 
{\sc D. Cassani,  B. Ruf, C. Tarsi,}
Optimal Sobolev type inequalities in Lorentz spaces. {\sl Potential Anal. }
 {\bf 39} (2013), no. 3, 265--285. 

\bibitem{CW} 
{\sc F.  Catrina, Z. Wang},
On the Caffarelli-Kohn-Nirenberg inequalities: sharp constants,
existence (and nonexistence),
{\sl Comm. Pure Appl. Math.} {\bf 54} (2001), no 2, 229--258.

\bibitem{Cham}
{\sc G. R. Chambers}, Proof of the Log-Convex Density Conjecture, 	
{\sl JEMS}, to appear. ArXiv:1311.4012v3 

\bibitem{ChiHo}
{\sc N. Chiba, T. Horiuchi}, On radial symmetry and its breaking in the 
Caffarelli-Kohn-Nirenberg inequalities for $p=1$. 
{\sl Math. J. Ibaraki Univ.} {\bf 47} (2015), 49--63. 

\bibitem{CP} 
{\sc E. Colorado, I. Peral}, Eigenvalues and bifurcation for elliptic equations 
with mixed Dirichlet-Neumann boundary conditions related 
to Caffarelli-Kohn-Nirenberg inequalities. 
{\sl Topological Methods in Nonlinear Analysis}, 
Journal of the Juliusz Schauder Center,
{\bf 23}, (2004), 239 --273.

\bibitem{C} 
{\sc  G. Csat\'o}, 
An isoperimetric problem with density and the Hardy Sobolev inequality in 
${\mathbb R}^2$, {\sl Differential Integral Equations} {\bf 28} (2015), no. 9-10, 971--988.

\bibitem{DDNT} 
{\sc  J. Dahlberg, A. Dubbs, E. Newkirk, H. Tran},
Isoperimetric regions in the plane with density $r^p$, 
{\sl New York J. Math.} {\bf 16} (2010), 31--51. 

\bibitem{diGiosia_etal}
{\sc L. Di Giosia, J. Habib, L. Kenigsberg, D. Pittman, W. Zhu},  
Balls Isoperimetric in $\mathbb{R}^n $ with Volume and Perimeter Densities $r^m$  and $r^k$
(2016), ArXiv:1610.05830v1.


\bibitem{DHHT} 
{\sc  A. Diaz, N. Harman, S. Howe, D. Thompson}
Isoperimetric problems in sectors with density. {\sl Adv. Geom. } 
{\bf 12} (2012), 589--619.

\bibitem{DolEstLoss}
{\sc J. Dolbeault, M. Esteban, M. Loss}, Rigidity versus symmetry 
breaking via nonlinear flows on cylinders and euclidean spaces, 
{\sl Invent. Math.} {\bf 206} (2016), no. 2, 397--440.

\bibitem{FleRi} 
{\sc W.H. Fleming, R. Rishel}, An integral formula for total gradient variation. 
{\sl Arch. Math. (Basel)} {\bf 11} (1960),  218--222. 

\bibitem{GaPe}
{\sc J.P. Garcia Azorero, I. Peral Alonso}, 
Hardy inequalities and some critical elliptic and
parabolic problems. {\sl J. of Differential Equations} {\bf 144} (1998), 444--476.

\bibitem{G} 
{\sc E. Giusti }, Minimal surfaces and functions of bounded variation. {\sl Monographs
in Mathematics, 80. Birkh\"{a}user Verlag}, Basel, 1984.



\bibitem{H} 
{\sc T. Horiuchi}, Best constant in weighted Sobolev inequality with weights
 being powers of distance from the origin. 
 {\sl J. Inequal. Appl.} {\bf 1} (1997), no. 3, 275--292.

\bibitem{HK} 
{\sc T. Horiuchi, P. Kumlin},
 On the Caffarelli-Kohn-Nirenberg-type inequalities 
involving critical and supercritical weights.
{\sl Kyoto J. Math.} {\bf 52} (2012), no. 4, 661--742.

\bibitem{Howe}
{\sc S. Howe}, The Log-Convex Density Conjecture and vertical surface 
area in warped products. 
{\sl Adv. Geom.} {\bf 15} (2015), 455--468.

\bibitem{Kaw}
{\sc B. Kawohl}, Rearrangements and convexity of level sets. Springer-Verlag N.Y. (1985). 

\bibitem{KZ} 
{\sc A.V. Kolesnikov, R.I. Zhdanov}, 
On isoperimetric sets of radially symmetric measures. Concentration, functional inequalities and isoperimetry, 123–-154,
{\sl Contemp. Math.} {\bf 545}, Amer. Math. Soc., Providence, RI, 2011.

\bibitem{Lan} 
{\sc R. Landes}, Some remarks on rearrangements and functionals with non-constant density. 
{\sl Math. Nachr.} {\bf 280} (2007), no. 5-6, 560--570.



\bibitem{MadernaSalsa}
{\sc C. Maderna, S. Salsa}, Sharp estimates for solutions to a
certain type of singular elliptic boundary value problems in two
dimensions.
{\sl Appl. Anal.} {\bf 12} (1981), 307--321.

\bibitem{M} 
{\sc V. Maz'ja}, Lectures on isoperimetric and isocapacitary inequalities
 in the theory of Sobolev spaces. Heat kernels and analysis on manifolds, 
 graphs, and metric spaces (Paris, 2002), 307--340, 
{\sl Contemp. Math.} {\bf 338}, Amer. Math. Soc., Providence, RI, 2003. 


\bibitem{Mo} 
{\sc F. Morgan}, Manifolds with density. {\sl Notices
Amer. Math. Soc.} \textbf{52} (2005), no.8, 853--858.

\bibitem{Mo1} 
{\sc F. Morgan}, Geometric Measure Theory: a Beginner's Guide. 
Academic Press, fifth edition, 2016.

\bibitem{Mo2} 
{\sc F. Morgan}, The Log-Convex Density Conjecture. {\sl
Contemporary Mathematics} \textbf{545} (2011), 209--211.

\bibitem{Musina} 
{\sc R. Musina}, Weighted Sobolev spaces of radially symmetric functions. 
{\sl Ann. Mat. Pura Appl.} (4) {\bf 193} (2014), no. 6, 1629--1659.

\bibitem{Stein} 
{\sc E. Stein}, 
Singular integrals and differentiability properties of functions. 
{\sl Princeton Mathematical Series}, no. 30, Princeton University Press, 
Princeton, N.J. 1970.

\bibitem{Talenti1} 
{\sc G. Talenti}, 
The standard isoperimetric theorem. {\sl Handbook of convex geometry}, 
Vol. A, B, 73--123, North-Holland, Amsterdam, 1993. 

\bibitem{Talenti2} 
{\sc G. Talenti},
A weighted version of a rearrangement inequality. Ann. Univ. Ferrara Sez. VII (N.S.) \textbf{43} (1997),  
(1998) 121--133.


\end{thebibliography}
\end{document}